\documentclass[11pt]{amsart}
\usepackage{amssymb, amstext, amscd, amsmath}
\usepackage[all]{xy}

\theoremstyle{plain}
\newtheorem{theorem}{Theorem}[section]
\newtheorem{corollary}[theorem]{Corollary}
\newtheorem{proposition}[theorem]{Proposition}
\newtheorem{lemma}[theorem]{Lemma}

\theoremstyle{definition}
\newtheorem{definition}[theorem]{Definition}
\newtheorem{example}[theorem]{Example}

\newtheorem{remark}[theorem]{Remark}

\makeatletter
\def\@cite#1#2{{\m@th\upshape\bfseries%
[{#1\if@tempswa{\m@th\upshape\mdseries, #2}\fi}]}} \makeatother

\newcommand{\bbC}{{\mathbb{C}}}
\newcommand{\bbN}{{\mathbb{N}}}
\newcommand{\bbZ}{{\mathbb{Z}}}

\newcommand{\B}{{\mathcal{B}}}
\newcommand{\C}{{\mathcal{C}}}

\newcommand{\F}{{\mathcal{F}}}
\newcommand{\G}{{\mathcal{G}}}

\newcommand{\J}{{\mathcal{J}}}
\newcommand{\K}{{\mathcal{K}}}
\renewcommand{\L}{{\mathcal{L}}}
\newcommand{\M}{{\mathcal{M}}}

\renewcommand{\O}{{\mathcal{O}}}

\newcommand{\T}{{\mathcal{T}}}

\newcommand{\X}{{\mathcal{X}}}

\newcommand{\fA}{{\mathfrak{A}}}
\newcommand{\fS}{{\mathfrak{S}}}

\renewcommand{\phi}{\varphi}
\newcommand{\upchi}{{\raise.35ex\hbox{\ensuremath{\chi}}}}
\def\gs{\sigma}
\def\ga{\alpha}
\def\gl{\lambda}

\newcommand{\foral}{\text{ for all }}

\newcommand{\alg}{\operatorname{alg}}
\newcommand{\id}{{\operatorname{id}}}
\newcommand\Span{\mathop{\rm span}}

\newcommand{\ca}{\mathrm{C}^*}
\newcommand{\cenv}{\mathrm{C}^*_{\text{env}}}
\newcommand{\sca}[1]{\left\langle#1\right\rangle} % \sca{a,b} =<a,b>
\newcommand{\lsca}[1]{\left[#1\right]}            % \lsca{a,b} =(a,b)
\newcommand{\nor}[1]{\left\Vert #1\right\Vert}    %\nor{x}=||x||

\newcommand{\vrt}{\G^{(0)}}
\newcommand{\vrtm}{\G^{(0)}_{-}}
\newcommand{\edg}{\G^{(1)}}
\newcommand{\xtau}{X_{\tau}}
\newcommand{\atau}{A_{\tau}}
\newcommand{\ptau}{\phi_{\tau}}

\newcommand{\sse}{\stackrel{\text{s}}{\thicksim}}
\newcommand{\sme}{\stackrel{\text{SME}}{\thicksim}}
\newcommand{\se}{\stackrel{\text{SE}}{\thicksim}}
\newcommand{\tsse}{{\stackrel{\text{SSE}}{\thicksim}}}

\begin{document}

%%%%%%%%%%%%%%%%%%%%%%%%%%%%%%%%%%%%%
\title[Operator algebras and C*-correspondences: A survey]{Operator algebras and C*-correspondences:\\ A survey}

\author[E.T.A. Kakariadis]{Evgenios T.A. Kakariadis}
\address{Pure Mathematics Department, University of Waterloo,
   Ontario N2L-3G1\\
   Canada}
\email{ekakaria@uwaterloo.ca}

\author[E.G. Katsoulis]{Elias~G.~Katsoulis}
\address{ Department of Mathematics\\University of Athens
\\ 15784 Athens \\GREECE \vspace{-2ex}}
\address{\textit{Alternate address:} Department of Mathematics
\\East Carolina University\\ Greenville, NC 27858\\USA}
\email{katsoulise@ecu.edu}


\maketitle
%%%%%%%%%%%%%%%%%%%%%%%%%%%%%%%%%%%%%

%%%%%%%%%%%%%%%%%%%%%%%%%%%%%%%%%%%%%
\begin{abstract}
In this paper we survey our recent work on $\ca$- correspondences and their associated operator algebras; in particular, on adding tails, the Shift Equivalence Problem and Hilbert bimodules.
\end{abstract}

%%%%%%%%%%%%%%%%%%%%%%%%%%%%%%%%%%%%%
\section*{Introduction}
%%%%%%%%%%%%%%%%%%%%%%%%%%%%%%%%%%%%%

Traditionally, the theory of $\ca$-correspondences has been used to generalize concrete results either from the theory of Cuntz-Krieger $\ca$-algebras or from the theory of crossed product $\ca$-algebras. Our goal in this project has been to discover a path in the opposite direction: to use the theory of $\ca$-correspondences in order to obtain results which are new even for the Cuntz-Krieger or the crossed product $\ca$-algebras. For instance in Theorem \ref{non-injective case} we describe a general process for adding  tail to a $\ca$-correspondence that has lead to new results in the theory of crossed product $\ca$-algebras and Cuntz's twisted crossed products. In the same spirit, our Theorem \ref{full} on the concept of shift equivalence of $\ca$-correspondences leads to a new result regarding the strong Morita equivalence of Cuntz-Krieger $\ca$-algebras (Theorem \ref{mor equiv}(4)). Therefore, even though the next few pages may seem rather abstract or encyclopedic, we believe that this abstraction will benefit even the reader who is interested only in the special classes of operator algebras mentioned above.

%%%%%%%%%%%%%%%%%%%%%%%%%%%%%%%%%%%%%
\section{Preliminaries}

Let $A$ be a $\ca$-algebra. An \emph{inner-product right $A$-module} is a linear space $X$ which is a right $A$-module together with a map
\begin{align*}
(\cdot,\cdot)\mapsto
\sca{\cdot,\cdot}_X \colon X \times X \rightarrow A
\end{align*}
such that
\begin{align*}
&\sca{\xi, \gl y + \eta}_X= \gl\sca{x,y}_X + \sca{x,\eta}_X \qquad
&(\xi,
y,\eta \in X, \gl\in \bbC),\\
&\sca{\xi, \eta a}_X= \sca{\xi,\eta}_X a \qquad &(\xi,\eta \in X,
a\in A),\\
&\sca{\eta,\xi}_X= \sca{\xi,\eta}_X^* \qquad &(\xi, \eta \in X),\\
&\sca{\xi,\xi}_X \geq 0 \text{; if }\sca{\xi,\xi}_X= 0 \text{ then }
\xi=0.
\end{align*}
A compatibility relation for the scalar multiplication is required, that is $\gl(\xi a)= (\gl \xi) a= \xi (\gl a)$, for all $\gl \in \bbC, \xi \in X, a\in A$. For $\xi \in X$ we write $\nor{\xi}_X:= \nor{\sca{\xi,\xi}_A}_A$ and one can deduce that $\nor{\cdot}_X$ is actually a norm. $X$ equipped with that norm will be called \emph{right Hilbert $A$-module} if it is complete and will be denoted as $X_A$.

For $X, Y$ Hilbert $A$-modules we define the set $\L(X,Y)$ of the \emph{adjointable maps} that consists of all maps $s:X \rightarrow Y$ for which there is a map $s^*: Y \rightarrow X$ such that
\begin{align*}
&\sca{s\xi,y}_Y= \sca{\xi,s^*y}_Y, & (x\in X, y\in Y).
\end{align*}
Every element of $\L(X,Y)$ is automatically a bounded $A$-linear map. An element $u \in \L(X,Y)$ is called \emph{unitary} if it is onto $Y$ and $\sca{u\xi, u\zeta}_Y= \sca{\xi,\zeta}_X$, for all $\xi, \zeta \in X$.

In particular, for $\xi\in X$ and $y\in Y$ we define $\Theta_{y,\xi}: X \rightarrow Y$  such that $\Theta_{y,\xi}(\zeta)= y \sca{\xi,\zeta}_X$, for all $\zeta \in X$. It is easy to check that $\Theta_{y,\xi} \in \L(X,Y)$ with $\Theta_{y,\xi}^*=\Theta_{\xi,y}$. We denote by $\K(X,Y)$ the closed linear space of $\L(X,Y)$ spanned by $\{\Theta_{y,\xi}: \xi\in X, y \in Y\}$. If $X=Y$ then $\K(X,X)\equiv \K(X)$ is a closed ideal of the $\ca$-algebra $\L(X,X)\equiv \L(X)$.

In a dual way we call $X$ a \emph{left Hilbert $A$-module} if it is complete with respect to the norm induced by an \emph{inner-product left $A$-module} $\lsca{\cdot,\cdot}_X$. The term \emph{Hilbert module} is reserved for the right Hilbert modules, whereas the left case will be clearly stated.

Given a Hilbert $A$-module $X$ over $A$, let $X^*= \{\xi^* \in \L(X,A): \xi^*(\zeta)= \sca{\xi,\zeta}_X \}$ be the \emph{dual left Hilbert $A$-module}, with
\begin{align*}
& a\cdot \xi^*= (\xi a^*)^*, & (\xi \in X, a\in A),\\
& \lsca{\xi^*,\zeta^*}_{X^*}= \sca{\xi,\zeta}_X, & (\xi, \zeta \in
X).
\end{align*}

A Hilbert $A$-module $X$ is called \emph{self-dual} when $X^*$ coincides with the set of bounded  (not necessarily adjointable) $A$-linear mappings from $X$ to $A$, hence a Riesz-Fr\'{e}chet Theorem is valid.

%%%%%%%%%%%%%%%%%%%%%%%%%%%%%%%%
\begin{example}
\textup{A $\ca$-algebra $A$ is a (trivial) Hilbert
$A$-module, when it is viewed as a Banach space with $a\cdot b:= ab$ and $\sca{a,b}_{A}:= a^*b$ for all $a,b\in A$.
It is a left Hilbert $A$-module when it is endowed with the left inner product $\lsca{a,b}_{A}:= ab^*$, for all $a,b\in A$.
Finally $\K(A) \simeq A$ and $\L(A_A) \simeq \M(A)$, i.e. the multiplier algebra of $A$.}
\end{example}

%%%%%%%%%%%%%%%%%%%%%%%%%%%%%%%%%%%%%
\subsection{$\ca$-correspondences}
%%%%%%%%%%%%%%%%%%%%%%%%%%%%%%%%%%%%%

Even though the class of $\ca$~-~correspondences has been thoroughly investigated for the last 25 years, the terminology still differs from author to author. We therefore present the terminology that we will be using in this paper.

%%%%%%%%%%%%%%%%%%%%%%%%%%%%%%%%%%%%%
\begin{definition}
\textup{An \emph{$A$-$B$-correspondence $X$} is a right Hilbert $B$-module together with a $*$-homomorphism $\phi_X\colon A \rightarrow \L(X)$. We will denote this by ${}_A X_B$. When $A=B$ we will simply refer to $X$ as \emph{$\ca$-correspondence over $A$}.}

\textup{A submodule $Y$ of $X$ is a \emph{subcorrespondence} of ${}_A X_B$, if it is a $C$-$D$- correspondence for some $\ca$-subalgebras $C$ and $D$ of $A$ and $B$, respectively.}
\end{definition}

A $\ca$-correspondence over $A$ is called \emph{non-degenerate} if the closed linear span of $\phi_X(A)X$ is dense in $X$. Moreover, $X$ is called \emph{full} if $\sca{X,X}_X$ is dense in $A$. Also, $X$ is called \emph{regular} if it is \emph{injective}, i.e. $\phi_X$ is injective, and $\phi_X(A) \subseteq \K(X)$.

Two $A$-$B$-correspondences $X$ and $Y$ are called unitarily equivalent, if there is a unitary $u \in \L(X,Y)$ such that $u(\phi_X(a)\xi b)= \phi_Y(a)(u(\xi))b$, for all $a\in A, b\in B, \xi \in X$. In that case we write $X \approx Y$. We write $X \lesssim Y$ when $X\approx Y_0$ for a subcorrespondence $Y_0$ of $Y$.

%%%%%%%%%%%%%%%%%%%%%%%%%%%%%%%%%%%%%
\begin{example}
\textup{Every Hilbert $A$-module $X$ is a $\K(X)$-$A$-correspondence when endowed with the left multiplication $\phi_X\equiv \id_{\K(X)}: \K(X) \rightarrow \L(X)$. A left inner product over $\K(X)$ can be defined by
$\lsca{\xi, \eta}_{X}= \Theta_{\xi, \eta}$, for all $\xi, \eta \in X$.
Also $X^*$ is an $A$-$\K(X)$-correspondence, when endowed with the following operations}
\begin{align*}
&\sca{\xi^*, \eta^*}_{X^*}=\lsca{\xi,\eta}_X= \Theta_{\xi, \eta} &(\xi^*, \eta^* \in X^*),\\
& \xi^* \cdot k =(k^*\xi)^*  &(\xi^*\in X^*, k\in \K(X)),\\
&\phi_{X^*}(a)\xi^*= a\cdot \xi^*= (\xi \cdot a^*)^* &(\xi^* \in X^*, a\in A).
\end{align*}
\end{example}

%%%%%%%%%%%%%%%%%%%%%%%%%%%%%%%%%%%%%
\begin{example}\label{E:K(X,Y)}
\textup{For Hilbert $A$-modules $X$ and $Y$, $\L(X,Y)$ becomes $\L(Y)$-$\L(X)$-correspondence by defining $\sca{s,t}:=s^*t$, $t\cdot a:= ta$ and $b\cdot t:=bt$, for every $s,t\in \L(X,Y), a\in \L(X)$ and $b\in \L(Y)$.}

\textup{Trivially, $\K(X,Y)$ is a $\K(Y)$-$\K(X)$-subcorrespondence of $\L(X,Y)$. Note that, when a c.a.i. in $\sca{X,X}_X$ is a right c.a.i. for $Y$, then $\K(Y)$ acts faithfully on $\K(X,Y)$. When $X=Y$ this is automatically true.}
\end{example}

%%%%%%%%%%%%%%%%%%%%%%%%%%%%%%%%%%%%%
\subsection{Interior Tensor Product}
%%%%%%%%%%%%%%%%%%%%%%%%%%%%%%%%%%%%%

The interior tensor product of two Hilbert modules plays the role of a generalized multiplication, stabilized by the elements of a common $\ca$-algebra (see \cite{Lan95} for more details). Let the $\ca$-correspondences ${}_A X_B$ and ${}_B Y_C$; the \emph{interior} or \emph{stabilized tensor product}, denoted by $X \otimes_B Y$ or simply by $X \otimes Y$, is the quotient of the vector space tensor product $X \otimes_{\alg} Y$ by the subspace generated by the elements of the form
\begin{align*}
\xi a \otimes y - \xi \otimes \phi(a)y, \foral \xi \in X, y \in Y, a
\in A.
\end{align*}
It becomes a Hilbert $B$-module when equipped with
\begin{align*}
& (\xi \otimes y)b:= \xi \otimes (y b), & (\xi \in X, y\in Y, b
\in B),\\
& \sca{\xi_1\otimes y_1, \xi_2\otimes y_2}_{X\otimes Y}:= \sca{y_1,
\phi(\sca{\xi_1,\xi_2}_X)y_2}_Y, & (\xi_1, \xi_2 \in X, y_1, y_2 \in
Y).
\end{align*}
For $s\in \L(X)$ we define $s\otimes \id_Y \in \L(X\otimes Y)$ as the mapping $\xi\otimes y \mapsto s(\xi)\otimes y$. Clearly, $(s\otimes \id_Y)^*= s^* \otimes \id_Y$; hence $X \otimes Y$ becomes an $A$-$C$-correspondence by defining $\phi_{X\otimes Y}(a):= \phi_X(a) \otimes \id_Y$. One can prove that the interior tensor product is associative, that is if $Z$ is a $C$-$D$-correspondence, then $(X \otimes_B Y) \otimes_C Z = X \otimes_B (Y \otimes_C Z)$.

%%%%%%%%%%%%%%%%%%%%%%%%%%%%%%%%%%%%%
\begin{example}\label{Ex:x^* otimes x}
\textup{When a Hilbert $A$-module $X$ is considered as the $\K(X)$-$A$- correspondence, then
$X \otimes_A X^* \approx \K(X)$ as $\ca$-correspondences over $\K(X)$, via the mapping}
\[
u_1: X \otimes_{A} X^* \rightarrow \K(X) : \xi \otimes \zeta^* \mapsto \Theta_{\xi,\zeta},
\]
\textup{and $X^* \otimes_{\K(X)} X \approx \sca{X,X}_X$, as $\ca$-correspondences over $A$ via the mapping}
\[
u_2: X^*\otimes_{\K(X)} X \rightarrow \sca{X,X}_X : \xi^*\otimes \zeta \mapsto \sca{\xi,\zeta}
\]
\textup{In particular $X^*\otimes_{\K(X)} X \approx A$, when $X$ is full.}
\end{example}

%%%%%%%%%%%%%%%%%%%%%%%%%%%%%%%%%%%%%
\subsection{Hilbert Bimodules}\label{subsection hilbert bimodules}
%%%%%%%%%%%%%%%%%%%%%%%%%%%%%%%%%%%%%

There are ${}_A X _B$ $\ca$-correspondences that are both left and right Hilbert modules. If a compatibility relation is satisfied between the two inner products then it is called \emph{Hilbert bimodule}.

%%%%%%%%%%%%%%%%%%%%%%%%%%%%%%%%%%%%%
\begin{definition}
A \emph{Hilbert $A$-$B$-bimodule} ${}_A X_B$ is a $\ca$-correspondence $X$ together with a left inner product $\lsca{\cdot,\cdot}_X\colon X \times X \rightarrow A$, which
satisfy:
\begin{align*}
& \lsca{\phi_X(a)\xi,\eta}_X= a\lsca{\xi,\eta}_X, & (\xi, \eta\in X, a\in A),\\
& \lsca{\xi,\eta}_X=\lsca{\eta,\xi}^*_X, &(\xi, \eta \in X)\\
& \lsca{\xi,\xi}_X\geq 0 \text{; if }\lsca{\xi,\xi}_X= 0 \text{
then
} \xi=0, \\
& \phi_X(\lsca{\xi,\eta}_X)\zeta=\xi\sca{\eta,\zeta}_X, & (\xi,
\eta,\zeta \in X).
\end{align*}
\end{definition}
The last equation implies that
$\phi_X(\lsca{\xi,\eta}_X)=\Theta^X_{\xi,\eta}$.
It is clear that Hilbert bimodules are a special case of $\ca$-correspondences. Let $I_X$ be the ideal,
\[
 I_X=\overline{\Span}\{\lsca{\xi,\eta}_X: \xi,\eta \in X\},
\]
in $A$. Using the very definitions, one can prove that $a\in \ker\phi_X$, if and only if $a\in I_X^\bot$. Hence, $\phi_X$ is $*$-injective, if and only if, the Hilbert $A$-bimodule $X$ is \emph{essential}, i.e. when the ideal $I_X$ is essential in $A$.

\begin{definition}
An \emph{$A$-$B$-imprimitivity bimodule} or \emph{equivalence bimodule} is an $A$-$B$-bimodule $M$ which is simultaneously a full left and a full right Hilbert $A$-module. That is $\lsca{M,M}_M$ is dense in $A$ and $\sca{M,M}_M$ is dense in $B$.
\end{definition}
It is easy to see that when $X$ is an $A$-$B$-imprimitivity bimodule then $A\simeq^{\phi} \K(M)$. Thus imprimitivity bimodules are automatically non-degenerate and regular.

%%%%%%%%%%%%%%%%%%%%%%%%%%%%%%%%%%%%%
\subsection{Matrix $\ca$-correspondences}
%%%%%%%%%%%%%%%%%%%%%%%%%%%%%%%%%%%%%

There is a number of ways of considering \emph{a} direct sum of $\ca$-correspondences. They are contained (as subcorrespondences) in the notion of the \emph{matrix $\ca$-correspondence} that is presented below. For the $\ca$-correspondences ${}_A E_A$, ${}_B F_B$, ${}_A R_B$ and ${}_B S_A$ the \emph{matrix $\ca$-correspondence}
$X= \left[\begin{array}{c|c} E & R\\
S & F \end{array}\right]$
over $A\oplus B$ is the Hilbert $(A\oplus B)$-module of the linear space of the matrices
$\left[\begin{array}{c|c} e & r \\
s & f \end{array}\right]$,
$e\in E, r\in R, s\in S, f\in F$, with {\small
\begin{align*}
&\left[\begin{array}{c|c} e & r \\ s & f
\end{array}\right] \cdot (a,b)= \left[\begin{array}{c|c} ea & rb \\ sa &
fb
\end{array}\right],\\
& \sca{\left[\begin{array}{c|c} e_1 & r_1 \\ s_1 & f_1
\end{array}\right], \left[\begin{array}{c|c} e_1 & r_1 \\ s_1 & f_1
\end{array}\right]}_X= \big(\sca{e_1,e_2}_E + \sca{s_1,s_2}_S,
\sca{r_1,r_2}_R + \sca{f_1,f_2}_F\big),
\end{align*}
} such that the $*$-homomorphism $\phi\colon A\oplus B \rightarrow \L\left(\left[\begin{array}{c|c} E & R\\ S & F \end{array}\right]\right)$ is defined as follows
\begin{align*}
\phi(a,b) \left[\begin{array}{c|c} e & r \\ s & f
\end{array}\right]= \left[\begin{array}{c|c} \phi_E(a)e & \phi_R(a)r
\\ \phi_S(b)s & \phi_F(b)f
\end{array}\right].
\end{align*}

The careful reader can see that this is exactly the \emph{exterior direct sum $\ca$-correspondence} of the two \emph{interior direct sum $\ca$-correspondences} ${}_{A\oplus B} (E+S)_A$ and ${}_{A\oplus B} (R+F)_B$. Hence, the linear space of the matrices
$\left[\begin{array}{c|c} e & r \\
s & f
\end{array}\right]$, with $e\in E, r\in R, s\in S, f\in F$, is complete with respect to the induced norm, thus a Hilbert $(A\oplus B)$-module. Moreover $E, F, R, S$ imbed naturally as subcorrespondences in $\left[\begin{array}{c|c} E & R\\ S & F \end{array}\right]$. The following lemma reasons the terminology \emph{matrix $\ca$-correspondence}, as tensoring comes by ``multiplying'' the matrices.

\begin{lemma}\label{technical lemma}
Let $E,F,R,S$ be $\ca$-correspondences as above. Then
\begin{align*}
&\left[\begin{array}{c|c} E & R\\ S & F \end{array}\right]
\otimes_{A\oplus B} \left[\begin{array}{c|c} E & R\\ S & F
\end{array}\right] \approx\\ &\approx \left[\begin{array}{c|c} (E\otimes_A E) +
(R\otimes_B S) & (E\otimes_A R) + (R \otimes_B F) \\
(S\otimes_A E) + (F\otimes_B S) & (S\otimes_A R) + (F\otimes_B F)
\end{array}\right],
\end{align*}
(unitary equivalent) as $\ca$-correspondences.
\end{lemma}

%%%%%%%%%%%%%%%%%%%%%%%%%%%%%%%%%%%%%
\subsection{Representations of $\ca$-correspondences}
%%%%%%%%%%%%%%%%%%%%%%%%%%%%%%%%%%%%%

Let us make a brief presentation on the representation theory of $\ca$-correspondences. Let $X$ be a $\ca$-correspondence over $A$. A (Toeplitz) representation $(\pi,t)$ of $X$ into a $\ca$-algebra $B$, is a pair of a $*$-homomorphism $\pi\colon A \rightarrow B$ and a linear map $t\colon X \rightarrow B$, such that
\begin{enumerate}
 \item $\pi(a)t(\xi)=t(\phi_X(a)(\xi))$,
 \item $t(\xi)^*t(\eta)=\pi(\sca{\xi,\eta}_X)$,
\end{enumerate}
for $a\in A$ and $\xi,\eta\in X$. An easy application of the $\ca$-identity shows that $t(\xi)\pi(a)=t(\xi a)$ is also valid. A representation $(\pi , t)$ is said to be \textit{injective} if $\pi$ is injective; in that case $t$ is an isometry.

The $\ca$-algebra generated by a representation $(\pi,t)$ equals the closed linear span of $t^n(\bar{\xi})t^m(\bar{\eta})^*$, where for simplicity we used the notation $\bar{\xi}\equiv \xi_{1}\otimes \dots \otimes \xi_{n} \in X^{\otimes n}$ and $t^n(\bar{\xi})\equiv t(\xi_1)\dots t(\xi_n)$. For any representation $(\pi,t)$ there exists a $*$-homomorphism $\psi_t:\K(X)\rightarrow B$, such that $\psi_t(\Theta^X_{\xi,\eta})= t(\xi)t(\eta)^*$.

Let $J$ be an ideal in $\phi_X^{-1}(\K(X))$; we say that a representation $(\pi,t)$ is $J$-coisometric if
\[
 \psi_t(\phi_X(a))=\pi(a), \text{ for any } a\in J.
\]
The representations $(\pi,t)$ that are $J_{X}$-coisometric, where
\[
 J_X=\ker\phi_X^\bot \cap \phi_X^{-1}(\K(X)),
\]
are called \emph{covariant representations}~\cite{Kats04}. The ideal $J_X$ is the largest ideal on which the restriction of $\phi_X$ is injective.

We define the Toeplitz-Cuntz-Pimsner algebra $\T_X$ as the universal $\ca$-algebra for all Toeplitz representations of $X$. Similarly, the Cuntz-Pimsner algebra $\O_X$ is the universal $\ca$-algebra for all covariant representations of $X$.

A concrete presentation of both $\T_{X}$ and $\O_{X}$ can be given in terms of the generalized Fock space $\F_{X}$ which we now describe. The \emph{Fock space} $\F_{X}$ over the correspondence $X$ is the interior direct sum of the $X^{\otimes n}$ with the structure of a direct sum of $\ca$-correspondences over $A$,
\[
\F_{X}= A \oplus X \oplus X^{\otimes 2} \oplus \dots .
\]
Given $\xi \in X$, the (left) creation operator $t_{\infty}(\xi) \in \L(\F_{X})$ is defined as
\[
t_{\infty}(\xi)(\zeta_0 , \zeta_{1}, \zeta_{2}, \dots ) = (0, \xi \zeta_0, \xi
\otimes \zeta_1, \xi \otimes \zeta_2, \dots),
\]
where $\zeta_n \in X^{\otimes n}$, $n \geq 0$ and $\zeta_0\in A$. (Here $X^{\otimes 0}\equiv A$, $X^{\otimes 1} \equiv X$ and $X^{\otimes n}= X \otimes X^{\otimes n-1}$, for $n\geq 2$.) For any $a \in A$, we define $\pi_{\infty}(a)\in \L(\F_{X})$ to be the diagonal operator with $\phi_X(c)\otimes \id_{n-1}$ at its $X^{\otimes n}$-th entry. It is easy to verify that $( \pi_{\infty}, t_{\infty})$ is a representation of $X$ which is called the \emph{Fock representation} of $X$. Fowler and Raeburn \cite{FR} (resp. Katsura \cite{Kats04}) have shown that the $\ca$-algebra $\ca( \pi_{\infty}, t_{\infty})$ (resp $\ca( \pi_{\infty}, t_{\infty})/ \K(\F_{X})J_{X}$) is $*$-isomorphic to $\T_{X}$ (resp. $\O_{X}$).

%%%%%%%%%%%%%%%%%%%%%%%%%%%%%%%%%%%%%
\begin{definition}
The \emph{tensor algebra} $\T_{X}^+$ of a $\ca$-correspondence $_{A} X_A$ is the norm-closed subalgebra of $\T_X$ generated by all elements of the form $\pi_{\infty}(a), t_{\infty}^n(\bar{\xi})$, $a \in A$, $\bar{\xi} \in X^n$, $n \in \bbN$.
\end{definition}

The tensor algebras for $\ca$-correspondences were pioneered by Muhly and Solel in \cite{MS}. They form a broad class of non-selfadjoint operator algebras which includes as special cases Peters' semicrossed products \cite{Pet2}, Popescu's non-commutative disc algebras \cite{Pop3}, the tensor algebras of graphs (introduced in \cite{MS} and further studied in \cite{KaK}) and the tensor algebras for multivariable dynamics \cite{DavKat}, to mention but a few. For more examples, see \cite{Kats03}.

There is an important connection between $\T_X^+$ and $\O_X$ given in the following Theorem of Katsoulis and Kribs \cite{KatsKribs06}. Recall that, for an operator algebra $\fA$ and a completely isometric representation $\iota \colon \fA \rightarrow A$, where $A=\ca(\iota(\fA))$, the pair $(A,\iota)$ is called \emph{a $\ca$-cover for $\fA$}. The \emph{$\ca$-envelope of the operator algebra $\fA$} is the universal $\ca$-cover $(A,\iota)$ such that, if $(B,\iota')$ is any other $\ca$-cover for $\fA$, then there exists a (unique) $*$-epimorphism $\Phi:B \rightarrow A$, such that $\Phi(\iota'(a))=\iota(a)$, for any $a\in \fA$ . For the existence of the $\ca$-envelope see \cite{Ham79, DrMc05, Ar06, BleLeM04, Kak11-2} .

\begin{theorem} \label{KKenv}
\textup{\cite[Theorem 3.7]{KatsKribs06}} The $\ca$-envelope of the tensor algebra $\T^+_X$ is $\O_X$.
\end{theorem}
As a consequence the Toeplitz-Cuntz-Pimsner algebra is the extension of the Cuntz-Pimsner algebra by the \v{S}ilov ideal. (Any ideal $\J \subseteq C$, for a $\ca$-cover $(C,\iota)$ of an operator algebra $\fA$, with the property that the restriction of the natural projection $C \rightarrow C/\J$ on $\fA$ is a complete isometry, is called a \textit{boundary ideal} and the \v{S}ilov ideal is the largest such ideal.)\\

Now, let us see how we can generalize the previous facts in the case of arbitrary $\ca$-correspondences by using the notion of matrix $\ca$-correspondences. A representation of an $A$-$B$-correspondence $X$ should be a triplet $(\pi_A, \pi_B, t)$ such that $\pi_A$ and $\pi_B$ are $*$-homomorphisms, $t$ is a linear mapping of $X$ and
\begin{enumerate}
 \item $\pi_A(a)t(\xi)\pi_B(b)=t(\phi_X(a)(\xi)b)$,
 \item $t(\xi)^*t(\eta)=\pi_B(\sca{\xi,\eta}_X)$,
\end{enumerate}
for all $\xi\in X, a\in A, b\in B$. If that is the case then one could define $\hat{t}:X \rightarrow \B(H^{(2)})$ such that
$\hat{t}(\xi)= \left[\begin{array}{cc} 0 & t(\xi)\\ 0 & 0
\end{array}\right]$. Then $(\pi_A \oplus \pi_B, \hat{t})$ defines a representation of the $(A\oplus B)$-correspondence
$\left[\begin{array}{c|c} 0 & X \\ 0 & 0 \end{array}\right]$.
Conversely, if $(\pi,t)$ is a representation of
$\left[\begin{array}{c|c} 0 & X \\ 0 & 0
\end{array}\right]$,
then $(\pi|_A, \pi|_B, t)$ defines a representation of $X$. Hence, we can identify $X$ with $\left[\begin{array}{c|c} 0 & X \\ 0 & 0 \end{array}\right]$ and define the Toeplitz-Cuntz-Pimsner, the Cuntz-Pimsner and the tensor algebra of the $A$-$B$-correspondence $X$ as the corresponding algebras of the $(A\oplus B)$-correspondence
$\left[\begin{array}{c|c} 0 & X \\ 0 & 0 \end{array}\right]$.  However, note that most of the results known for $\ca$-correspondences over the same $\ca$-algebra must be verified, basically because $X^{\otimes n}$ is absurd for all $n\geq 2$.

%%%%%%%%%%%%%%%%%%%%%%%%%%%%%%%%%%%%%
\begin{remark}
\textup{We already gave a brief description of $X^*$ of the Hilbert $A$-module $X$. When $X$ is a correspondence over $A$ this can be simplified. Let $(\pi_u,t_u)$ be the universal representation of ${}_A X_A$; then $X^*$ is the closed linear span of $t(\xi)^*, \xi \in X$ with the left multiplication and inner product inherited by the trivial correspondence $\ca(\pi_u,t_u)$. Via this identification one can produce a theory for the left analogue of $\ca$-correspondences (that means also, to start with left Hilbert modules), but in most of the cases it can be recovered.}
\end{remark}

%%%%%%%%%%%%%%%%%%%%%%%%%%%%%%%%%%%%%
\subsection{Examples}
%%%%%%%%%%%%%%%%%%%%%%%%%%%%%%%%%%%%%

One of the fundamental examples in the theory of $\ca$-correspondences  are the $\ca$-algebras of directed graphs. (For more details see~\cite{Raeb}.)

Let $\G$ be a countable directed graph with vertex set $\G^{(0)}$, edge set $\G^{(1)}$ and range and source maps $r$ and $s$ respectively. A family of partial isometries $\{ L_e\}_{ e\in \G^{(1)}}$ and projections $\{ L_p \}_{p\in \G^{(0)}}$ is said to obey the Cuntz-Krieger relations associated with $\G$ if and only if they satisfy
\[
(\dagger)  \left\{
\begin{array}{lll}
(1)  & L_p L_q = 0 & \mbox{$\forall\, p,q \in \G^{(0)}$, $ p \neq q$}  \\
(2) & L_{e}^{*}L_f = 0 & \mbox{$\forall\, e, f \in \G^{(1)}$, $e \neq f $}  \\
(3) & L_{e}^{*}L_e = L_{s(e)} & \mbox{$\forall\, e \in \G^{(1)}$}      \\
(4)  & L_e L_{e}^{*} \leq L_{r(e)} & \mbox{$\forall\, e \in \G^{(1)}$} \\
(5)  & \sum_{r(e)=p}\, L_e L_{e}^{*} = L_{p} & \mbox{$\forall\, p
\in \G^{(0)}$ with $|r^{-1}(p)|\neq 0 , \infty$}
\end{array}
\right.
\]
The relations $(\dagger)$ have been refined in a series of papers by the Australian school and reached the above form in \cite{BHRS, RS}. All refinements involved condition $(5)$ and as it stands now, condition $(5)$ gives the equality requirement for projections $L_{p}$ such that $p$ is not a source and receives finitely many edges. (Indeed, otherwise condition $(5)$ would not be a $\ca$-condition.)

It can been shown that there exists a universal $\ca$-algebra, denoted as $\O_{\G}$, associated with the relations $(\dagger)$. Indeed, one constructs a single family of partial isometries and projections obeying $(\dagger)$. Then, $\O_{\G}$ is the $\ca$-algebra generated by a `maximal' direct sum of such families. It turns out that there $\O_{\G}$ is the Cuntz-Pimsner algebra of a certain $\ca$-correspondence~\cite{MS}. The associated Cuntz-Pimsner -Toeplitz algebra is the universal algebra for the first four relations in $(\dagger)$ and is denoted as $\T_{\G}$.

%%%%%%%%%%%%%%%%%%%%%%%%%%%%%%%%%%%%%
\begin{example}
\textup{Let $G, G'$ be two graphs with adjacent matrices $A_G$ and $A_{G'}$. If $X_G$ and $X_{G'}$ are the corresponding $\ca$-correspondences, then $X_G \otimes X_{G'}$ is unitarily equivalent to the $\ca$-correspondence that comes from the adjacent matrix $A_G \cdot A_{G'}$}.

\textup{Let $X$ be an imprimitivity bimodule that comes from a graph $\G=(\vrt,\edg,r,s)$ (we follow the notation in \cite{Raeb}); then $\G$ is either a cycle or a double infinite path. Pick your favorite completely isometric representation of $\T_X^+$; for us it is a Cuntz-Krieger family $\{P_v, S_e: v\in \vrt, e \in \edg\}$ (because of \cite{KaK}). Apart from the usual relations, due to the form of the graph we have the simplified relation $S_eS_e^*=P_{r(e)}$, for all $e\in \edg$. Therefore $\{P_v, S_e^*: v\in \vrt, e\in \edg\}$ defines a Cuntz-Krieger family of the graph $\G^*=(\vrt,\edg,r^*,s^*)$ where the arrows are reversed, i.e. $r^*=s$ and $s^*=r$, thus $X^*$ is the $\ca$-correspondence coming from this graph $\G^*$.}
\end{example}

A second example comes from the class of dynamical systems. Let $\beta\colon A \rightarrow B$ be a $*$-homomorphism of $\ca$-algebras. The trivial Hilbert module $B_B$ becomes a $A$-$B$-correspondence, denoted by $(X_\beta,A)$, when endowed with the left action $\phi_B$ such that $\phi_B(a)b=\beta(a)b$ for all $a \in A$ and $b\in B$.

%%%%%%%%%%%%%%%%%%%%%%%%%%%%%%%%%%%%%
\begin{example}\label{Ex:ten dyn}
If $(X_\ga,C)$ is a $C$-$A$-correspondence via a $*$-homomorphism $\ga\colon C \rightarrow A$, and $(X_\beta,A)$ is a $A$-$B$-correspondence via a $*$-homomorphism $\beta\colon A \rightarrow B$, then $A\otimes_A B$ is unitarily equivalent to the $C$-$B$-correspondence $(X_{\beta\circ \ga},C)$ associated to the $*$-homomorphism $\beta \circ \ga\colon C \rightarrow B$.

Moreover, $(X_\beta,A)$ is an imprimitivity bimodule if and only if $\beta$ is a $*$-isomorphism. In this case $(X_\beta,A)^* \approx (X_{\ga^{-1}}, B)$. In particular, when $A=B$, then $\O_{X_\beta}$ is the usual crossed product $B\rtimes_\beta \bbZ$ and the tensor algebra $\T_{X_\beta}$ is Peter's semicrossed product \cite{Pet2} (for this and various types of semicrossed products see also \cite{kakar, Kak11, kakarkats10}).
\end{example}

%%%%%%%%%%%%%%%%%%%%%%%%%%%%%%%%%%%%%
\section{Adding Tails}\label{adding tails}
%%%%%%%%%%%%%%%%%%%%%%%%%%%%%%%%%%%%%

In \cite{KakKats11} the authors developed a method of ``adding tails'' that extends the one developed by Muhly and Tomforde in \cite{MuTom04}.\\

Let $\G$ be a connected, directed graph with a distinguished sink $p_0 \in \vrt$ and no sources. We assume that $\G$ is \emph{contractible} at $p_0$, i.e. the subalgebra $\bbC L_{p_0}$ is a full corner of the Cuntz-Krieger algebra $\O_\G$. So there exists a unique infinite path $w_0 = e_1 e_2 e_3\dots$ ending at $p_0$, i.e. $r(w_0)=p_0$. Let $p_n\equiv s(e_n)$, $n \geq 1$.

Let $(A_p)_{p \in \vrt}$ be a family of $\ca$-algebras parameterized by the vertices of $\G$ so that $A_{p_{0}}=A$. For each $e \in \edg$, we now consider a full, right Hilbert $A_{s(e)}$~-~module $X_e$ and a $*$-homomorphism
\[
\phi_e \colon A_{r(e)} \longrightarrow \L(X_e)
\]
satisfying the following requirements.

For $e \neq e_1$, the homomorphism $\phi_e$ are required to be injective and map onto $\K(X_e)$, i.e. $\phi_e(A_{r(e)})= \K(X_e)$. Therefore, each $X_e$, $e \neq e_1$, is an $A_{r(e)}$-$A_{s(e)}$-equivalence bimodule, in the language of Rieffel.

For $e = e_1$, we require $\K(X_{e_1}) \subseteq \phi_{e_1} (A)$ and
\begin{equation} \label{inj}
J_{X}\subseteq \ker \phi_{e_1} \subseteq \left( \ker \phi_X
\right)^{\perp}.
\end{equation}
In addition, there is also a \textit{linking condition}
\begin{equation} \label{linkin}
\phi_{e_1}^{-1}(\K(X_{e_1})) \subseteq \phi_X^{-1}(\K(X))
\end{equation}
required between the maps $\phi_{X}$ and $\phi_{e_1}$.

Let $T_0= c_0 (\, (A_p)_{p \in \vrtm})$ denote the $c_0$-sum of the family $(A_p)_{p \in \vrtm}$, where $\vrtm \equiv \vrt \backslash \{p_0\}$. Consider the set $c_{00}((X_e)_{e \in \edg}) \subseteq c_0 ((X_e)_{e \in \edg})$, consisting of sequences which are zero everywhere but on a finite number of entries. Equip $c_{00}((X_e)_{e \in \edg}) $ with a $T_0$-valued inner product defined by
\[
\sca{ u, v  }(p)= \sum_{s(e)=p} \, u^*_ev_e, \quad p \in \vrtm,
\]
for any $u, v \in c_{00}((X_e)_{e \in \edg})$. Let $T_1$ be the completion of $c_{00}((X_e)_{e \in \edg}) $ with respect to that inner product. Equip now $T_1$ with a right $T_0$~-~action, so that
\[
(ux)_e = u_ex_{s(e)}, \quad e \in \edg,
\]
for any $x \in T_0$, so that the pair $(T_1, T_0)$ becomes a right $T_0$-Hilbert module. The pair $(T_0, T_1)$ is the \textit{tail} for $_{A} X_A$.

To the $\ca$-correspondence $_{A} X_A$ and the data
\[
\tau\equiv \Big(\G, (X_e)_{e \in \edg}, (A_p)_{p\in \vrt},
(\phi_{e})_{e \in \edg} \Big),
\]
we now associate
\begin{align} \label{tau}
\begin{split}
\atau &\equiv A\oplus T_0  \\
\xtau &\equiv X \oplus  T_1.
\end{split}
\end{align}
Using the above, we view $\xtau$ as a $\atau$-Hilbert module with the standard right action and inner product for direct sums of Hilbert modules. We also define a left $\atau$-action $\ptau: \atau \rightarrow \L(\xtau)$ on $\xtau$ by setting
\[
\ptau(a, x \, )(\xi, u )= (\phi_X(a)\xi, v  ),
\]
where
\[
v_{e}= \left\{ \begin{array}{ll}
   \phi_{e_1}(a)(u_{e_1}),& \mbox{if $e = e_1$} \\
   \phi_e (x_{r(e)})u_e, & \mbox{otherwise}
   \end{array}
\right.
\]
for $a \in A$, $\xi \in X$, $x \in T_0$ and $u \in T_1$.

%%%%%%%%%%%%%%%%%%%%%%%%%%%%%%%%%%%%%
\begin{theorem} \label{non-injective case}
\cite[Theorem 3.10]{KakKats11} Let $_{A} X_A$ be a non-injective $\ca$- correspondence and let $\xtau$ be the graph $\ca$-correspondence over $\atau$ defined above. Then $\xtau$ is an injective $\ca$-correspondence and the Cuntz-Pimsner algebra $\O_X$ is a full corner of $\O_{\xtau}$.
\end{theorem}

Regarding Theorem~\ref{non-injective case} and the conditions imposed on the graph $\G$ and the maps $(\phi_e)_{e \in \edg}$, we have asked that the graph $\G$ be contractible. We cannot weaken this assumption to include more general graphs. Indeed, we want the tail associated with the data
 \[
\tau = \Big(\G, (X_e)_{e \in \edg}, (A_p)_{p\in \vrt}, (\phi_{e})_{e \in \edg} \Big),
 \]
to work with any possible Cuntz-Pimsner algebra $\O_X$ that can be ``added on''. This should apply in particular to the Cuntz-Krieger algebra $\O_{G_{p_0}}$ of the (trivial) graph $\G_{p_0}$ consisting only of one vertex $p_0$. By taking $\tau$ to be the ``usual'' tail associated with $\G$, i.e. $X_e=A_e=\bbC L_{p_0}$ and $\phi_e$ left multiplication for all $e$, we see that $\O_{G_{p_0}}$ is a full corner of $O_{\xtau}$ if and only if $\G$ is contractible at $p_0$.

Conditions~(\ref{inj}) and (\ref{linkin}) are also necessary, as the following result suggests.

%%%%%%%%%%%%%%%%%%%%%%%%%%%%%%%%%%%%%
\begin{proposition}
Let $_{A} X_A$ be a non-injective $\ca$-correspondence and let $\xtau$ be the $\ca$-correspondence over $\atau$ associated with the data
 \[
\tau = \Big(\G, (X_e)_{e \in \edg}, (A_p)_{p\in \vrt}, (\phi_{e})_{e
\in \edg} \Big),
 \]
as defined at the beginning of the section. If $X_\tau$ is injective,
\[
(\phi_1^{-1}(\K(X_1))+J_{X})\oplus 0 \subseteq
J_{\xtau},
\]
and the covariant representations of $X_\tau$ restrict to covariant representations of $\phi_X$, then
\begin{equation*}
 J_{X}\subseteq \ker \phi_{1}
\subseteq \left( \ker \phi_X \right)^{\perp},
\end{equation*}
and the \textit{linking condition}
\begin{equation*}
\phi_{1}^{-1}(\K(X_{1})) \subseteq \phi_X^{-1}(\K(X))
\end{equation*}
holds.
\end{proposition}

Lets see now how the work of Muhly and Tomforde fits in our theory.

%%%%%%%%%%%%%%%%%%%%%%%%%%%%%%%%%%%%%
\begin{example}[The Muhly-Tomforde tail~\cite{MuTom04}]
\textup{Given a (non-injective) correspondence $(X, A, \phi_X)$, Muhly and Tomforde construct in \cite{MuTom04} the tail that results from the previous construction, with respect to data
\[
\tau =\Big(\G, (X_e)_{e \in \edg}, (A_p)_{p\in \vrt}, (\phi_{e})_{e
\in \edg} \Big)
\]
defined as follows. The graph $\G$ is illustrated in the figure below.}

\vspace{.2in}
\[
\xymatrix{&{\bullet^{p_0}} &{\,\, \bullet}^{p_1}
\ar[l]^{e_1}&{\,\,\bullet^{p_2}} \ar[l]^{e_2} &{\bullet^{p_3}}
\ar[l]^{e_3}&{\,\,\bullet} \ar[l]&\dots \ar[l]}
\]
\vspace{.2in}

\noindent \textup{$A_p=X_e=\ker \phi_X$, for all $p \in \vrtm$ and $e \in \edg$. Finally,
\[
\phi_e(a)u_e=au_e,
\]
for all $e \in \edg$, $u_e \in X_e$ and $a \in A_{r(e)}$.}
\end{example}

%%%%%%%%%%%%%%%%%%%%%%%%%%%%%%%%%%%%%
\subsection{Applications}
%%%%%%%%%%%%%%%%%%%%%%%%%%%%%%%%%%%%%

%%%%%%%%%%%%%%%%%%%%%%%%%%%%%%%%%%%%%
\subsubsection{Semicrossed Products}
%%%%%%%%%%%%%%%%%%%%%%%%%%%%%%%%%%%%%

The tail of Muhly and Tomforde has had significant applications in the theory of $\ca$-correspondences, including a characterization for the $\ca$-envelope of the tensor algebra of a non-injective correspondence~\cite{KatsKribs06}. However, it also has its limitations, as we are about to see.

%%%%%%%%%%%%%%%%%%%%%%%%%%%%%%%%%%%%%
\begin{example}\label{MTexample}
Let $(X_{\ga}, A)$ be the $\ca$-correspondence canonically associated with a dynamical system $(A, \ga)$ and let $\O_{(A,\ga)}$ be the associated Cuntz-Pimsner $\ca$-algebra. If $\ga$ is not injective, then by using the Muhly-Tomforde tail we obtain an injective $\ca$-correspondence $(Y,B, \phi_Y)$ so that $\O_{(A,\ga)}$ is a full corner of $\O_Y$. Remarkably, $(Y, B, \phi_Y)$ may not come from any dynamical system in general. Assume that $\ker \ga \subseteq A$ is an essential ideal of $A$; then the Muhly-Tomforde tail produces an injective correspondence $(Y, B, \phi_Y)$ with
\begin{align*}
Y= A\oplus c_{0}(\ker \ga), \qquad B=A\oplus c_{0}(\ker \ga)
\end{align*}
and $\phi_Y$ defined by
\[
\phi_Y\big( a, (c_i)_i\big)\big(a',(c_i')_i\big)
= \big(\ga(a)a',\ga(a)c_1', c_1c_2', c_2c_3',\dots\big),
\]
where $a,a' \in A$ and $(c_i)_i, (c_i')_i \in c_0(\ker \ga)$.

If there was a $*$-homomorphism $\beta$ satisfying
\begin{equation} \label{impos}
\phi_Y(b)(b')=\beta(b)b',
\end{equation}
then by equating second coordinates in the equation
\[
\phi_Y\big( 1, (c_i)_i\big)\big(a',(c_i')_i\big)=
\beta \big( 1, (c_i)_i\big)\big(a',(c_i')_i\big)
\]
we would obtain,
\[
c_1'=
\beta \big(1,(c_i)_i \big)_2 c_1',
\]
for all $c_1' \in \ker \ga$. Since $\ker \ga$ is an essential ideal, we have $\ker \ga \ni \beta \big(1,(c_i)_i \big)_2$ $=1$, a contradiction.
\end{example}

Therefore, the Muhly-Tomforde tail produces an injective correspondence but not necessarily an injective dynamical system. Nevertheless, there exists a tail that can be added to $(X_{\ga}, A)$ and produce an injective correspondence that comes from a dynamical system.

%%%%%%%%%%%%%%%%%%%%%%%%%%%%%%%%%%%%%
\begin{example}
\textup{If $(X_{\ga}, A)$ is the $\ca$-correspondence canonically associated with a dynamical system $(A, \ga)$, then the appropriate tail for $(X_{\ga}, A)$ comes from the data
\[
\tau =\Big(\G, (X_e)_{e \in \edg}, (A_p)_{p\in \vrt}, (\phi_{e})_{e \in \edg} \Big)
\]
where $\G$ is as in Example~\ref{MTexample} \textit{but} for any $p \in \vrtm$ and $e \in \edg$,
\[A_p=X_e= \theta(A),
\]
where $\theta \colon A \longrightarrow M(\ker \ga)$ is the map that extends the natural inclusion $\ker \ga \subseteq M(\ker \ga))$ in the multiplier algebra. Finally
\[
\phi_e(a)u_e=au_e,
\]
for all $e \in \edg$, $u_e \in X_e$ and $a \in A_{r(e)}$. Then the correspondence $X_\tau$ is canonically associated to the dynamical system $(B,\beta)$, where}
\begin{align*}
B= A \oplus c_0(\theta(A) \text{ and } \beta(a,(c_i)_i)= (\ga(a),\theta(a), (c_i)_i).
\end{align*}
\end{example}

%%%%%%%%%%%%%%%%%%%%%%%%%%%%%%%%%%%%%
\subsubsection{Multivariable Dynamical Systems}
%%%%%%%%%%%%%%%%%%%%%%%%%%%%%%%%%%%%%

We now apply the method of adding tails to $\ca$-dynamics. Apart from their own merit, this application will also address the necessity of using more elaborate tails than that of Muhly and Tomforde in the process of adding tails to $\ca$-correspondences. This necessity has been already noted in the one-variable case.

A \textit{multivariable $\ca$-dynamical system} is a pair $( A, \ga)$ consisting of
a $\ca$-algebra $A$ along with a tuple $\ga=( \ga_1, \ga_2 , \dots, \ga_n)$, $n \in \bbN$, of $*$- endomorphisms of $A$. The dynamical system is called injective if $\cap_{i=1}^n \, \ker\ga_i =\{0\}$.

In the $\ca$-algebra literature, the algebras $\O_{(A, \ga)}$ are denoted as $A \times_{\ga} \O_n$ and go by the name "twisted tensor products by $\O_n$". They were first introduced and studied by Cuntz~\cite{Cun} in 1981. In the non-selfadjoint literature, these algebras are much more recent. In \cite{DavKat} Davidson and the second named author introduced the tensor algebra $\T_{(A, \ga)}$ for a  multivariable dynamical system $(A,\ga)$. It turns out that $\T_{(A,\ga)}$ is completely isometrically isomorphic to the tensor algebra for the $\ca$-correspondence $(X_{\ga}, A, \phi_{\ga})$. As such, $\O_{(A, \ga)}$ is the $\ca$-envelope of $\T_{(A, \ga)}$. Therefore, $\O_{(A, \ga)}$ provides a very important invariant for the study of isomorphisms between the tensor algebras $\T_{(A, \ga)}$.

To the multivariable system $( A, \ga)$ we associate a $\ca$-correspondence $(X_{\ga}, A, \phi_{\ga})$ as follows. Let $X_{\ga}=A^n = \oplus_{i=1}^n A$ be the usual right $A$-module. That is
\begin{enumerate}
\item $(a_1,\dots,a_n)\cdot a= (a_1 a ,\dots,a_n a)$,
\item $\sca{(a_1,\dots,a_n), (b_1,\dots,b_n)}=\sum_{i=1}^n
\sca{a_i,b_i}=\sum_{i=1}^n a_i^*b_i$.
\end{enumerate}
Also, by defining the $*$-homomorphism
\begin{align*}
\phi_{\ga}\colon A \longrightarrow \L(X_{\ga})\colon a \longmapsto \oplus_{i=1}^n \ga_i(a),
\end{align*}
$X_\ga$ becomes a $\ca$-correspondence over $A$, with $\ker\phi_{\ga}= \cap_{i=1}^n \ker \ga_i$ and $\phi(A) \subseteq \K(X_{\ga})$. It is easy to check that in the case where $A$ and all $\ga_i$ are unital, $X_\ga$ is finitely generated as an  $A$-module by the elements
\[
e_1:=(1,0,\dots,0),e_2:=(0,1,\dots,0),\dots,e_n:=(0,0,\dots,1),
\]
where $1\equiv 1_A$. In that case, $(\pi,t)$ is a representation of this $\ca$- correspondence if and only if $t(\xi_i)$ are isometries with pairwise orthogonal ranges and
\[
\pi(c)t(\xi)=t(\xi)\pi(\ga_i(c)), \quad i=1,\dots, n.
\]

%%%%%%%%%%%%%%%%%%%%%%%%%%%%%%%%%%%%%
\begin{definition}
The Cuntz-Pimsner algebra $\O_{(A, \ga)}$ of a multivariable $\ca$-dynamical system $(A,\ga)$ is the Cuntz-Pimsner algebra of the $\ca$- correspondence $(X_{\ga}, A, \phi_{\ga})$ constructed as above.
\end{definition}

The graph $\G$ that we associate with $(X_{\ga}, A, \phi_{\ga})$ has no loop edges and a single sink $p_0$. All vertices in $\vrt \backslash \{p_0\}$ emit $n$ edges, i.e. as many as the maps involved in the multivariable system, and receive exactly one. In the case where $n=2$, the following figure illustrates $\G$.
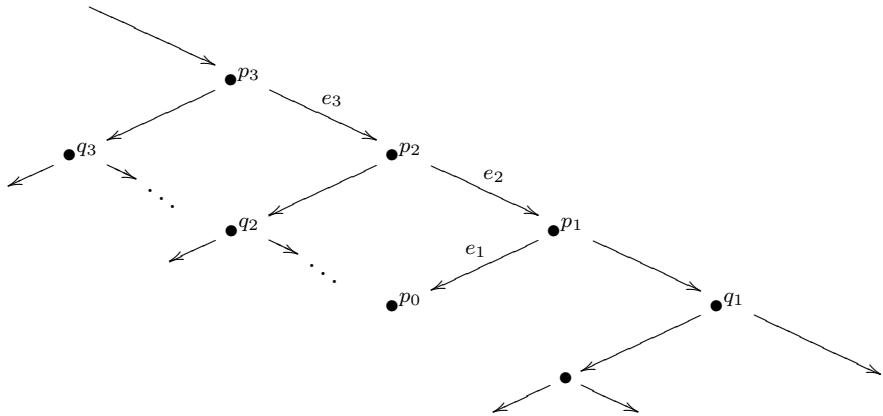
\vspace{.2in}
\[
\xygraph{
!{<0cm,0cm>;<2.15cm,0cm>:<0cm,1cm>::}
!{(1,3)}*+{\empty}="a"
!{(2,2)}*+{\bullet^{p_3}}="b"
!{(1,1)}*+{\bullet^{q_3}}="c"
!{(3,1)}*+{\bullet^{p_2}}="cc"
!{(0.5,.5)}*+{\empty}="e"
!{(1.5,.5)}*+{\ddots}="f"
!{(4,0)}*+{\bullet^{p_1}}="j"
!{(5,-1)}*+{\bullet^{q_1}}="q"
!{(6,-2)}*+{\empty}="u"
!{(7,-3)}*+{\empty}="bb"
!{(3,-1)}*+{\bullet^{p_0}}="p"
!{(4,-2)}*+{\bullet}="w"
!{(3.5,-2.5)}*+{\empty}="v"
!{(4.5,-2.5)}*+{\empty}="y"
!{(2,0)}*+{\bullet^{q_2}}="i"
!{(1.5,-.5)}*+{\empty}="k"
!{(2.5,-.5)}*+{\ddots}="l"
"b" :"cc"^{e_3}
"cc":"j"^{e_2}
"j":"q"
"j":"p"_{e_1}
"a":"b"
"q":"w"
"q":"u"
"cc":"i"
"b":"c"
"c":"e" "c":"f"
"i":"k"
"i":"l"
"w":"v"
"w":"y"
}
\]
Clearly, $\G$ is $p_0$-accessible. There is also a unique infinite path $w$ ending at $p_0$ and so $\G$ is contractible at $p_0$.\\

Let $\J \equiv \cap_{i=1}^{n} \, \ker \alpha_i$ and let $M(\J)$ be the multiplier algebra of $\J$. Let $\theta \colon A \longrightarrow M(\J)$ the map that extends the natural inclusion $\J \subseteq M(\J))$. Let $X_e= A_{s(e)} =\theta(A)$, for all $e \in \edg$, and consider $(X_e, A_{s(e)})$ with the natural structure that makes it into a right Hilbert module.

For $e \in \edg \backslash \{e_1\}$ we define $\phi_{e}(a)$ as left multiplication by $a$. With that left action, clearly $X_e$ becomes an $A_{r(e)}$-$A_{s(e)}$-equivalence bimodule. For $e = e_1$, it is easy to see that
\[
\phi_{e_1}(a)(\theta(b))\equiv \theta(ab), \quad a,b \in A
\]
defines a left action on $X_{e_1}=\theta(A)$, which satisfies both (\ref{inj}) and (\ref{linkin}).

%%%%%%%%%%%%%%%%%%%%%%%%%%%%%%%%%%%%%
\begin{theorem} \label{multitail}
\cite[Theorem 4.2]{KakKats11} If $(A, \ga)$ is a non-injective multivariable $\ca$-dynamical system, then there exists an injective multivariable $\ca$-dynamical system $(B, \beta)$ so that the associated Cuntz-Pimsner algebras $\O_{(A, \ga)}$ is a full corner of
$\O_{(B, \beta)}$. Moreover, if $A$ belongs to a class
$\, \C$ of $\, \ca$-algebras which is invariant under quotients and $c_0$-sums, then $B\in \C$ as well. Furthermore, if $(A, \ga)$ is non-degenerate, then so is $(B, \beta)$.
\end{theorem}

The reader familiar with the work of Davidson and Roydor may have noticed that the arguments in the proof of Theorem~\ref{multitail}, when applied to multivariable systems over commutative $\ca$-algebras produce a tail which is different from that of Davidson and Roydor in \cite[Theorem 4.1]{DavR}. It turns out that the proof of \cite[Theorem 4.1]{DavR} contains an error and the technique of Davidson and Roydor does not produce a full corner, as claimed in \cite{DavR}.  Nevertheless, \cite[Theorem 4.1]{DavR} is valid as Theorem \ref{multitail} demonstrates (a fact also mentioned in \cite[Corrigendum]{DavR}). We illustrate this by examining their arguments in the following simple case.

%%%%%%%%%%%%%%%%%%%%%%%%%%%%%%%%%%%%%
\begin{example}
\textup{ (The Davidson-Roydor tail~\cite{DavR}). Let $\X\equiv\{u ,v \}$ and consider the maps $\sigma_i \colon  \X \rightarrow \X$, $i=1,2$, with $\sigma_i(u) =v$ and $\sigma_i(v)=v$. Set $\sigma \equiv (\sigma_1 , \sigma_2)$ and let $\O_{(\X, \sigma)}$ be the Cuntz-Pimsner algebra associated with the multivariable system $(\X, \sigma)$, which by \cite{DavKat} is the $\ca$-envelope of the associate tensor algebra.}

\textup{We now follow the arguments of \cite{DavR}. In order to obtain $\O_{(\X, \sigma)}$ as a full corner of an injective Cuntz-Pimsner algebra, Davidson and Roydor add a tail to the multivariable system. They define
\[
T= \{(u,k)\mid k<0\} \mbox{ and }\X^{T}= \X\cup T.
\]
For each $1\leq i \leq2$, they extend $\sigma_i$ to a map $\gs_i^T\colon \X^T \rightarrow \X^T$ by
\[\gs^T(u,k)= (u, k+1) \mbox{ for $k<-1$ and  } \sigma_i^T(u,-1)=u.
\]
They then consider the new multivariable system $(\X^T, \sigma^T)$ and its associated Cuntz-Pimsner algebra $\O{(\X^T,\gs^T)}$.}

\textup{It is easy to see that the Cuntz-Pimsner algebra $\O_{(\X, \sigma)}$ for the multivariable system $(\X, \sigma)$ is the Cuntz-Krieger algebra $\O_{\G}$ of the graph $\G$ illustrated below, while the Cuntz-Pimsner algebra $\O{(\X^T,\gs^T)}$ is isomorphic to the Cuntz-Krieger algebra $\O_{\G^T}$ of the following graph $\G^{T}$, where for simplicity we write $u_k$ instead of $(u,k)$, $k<0$,}

\vspace{.2in}
{\small
\begin{align*}
\xymatrix{ \G: \hspace{-.7cm} &
 {\bullet^v} \ar@(u,ul)[] \ar@(d,dl)[] &{\, \bullet^u \,} \ar@/^/[l] \ar@/_/[l] }, \quad
\xymatrix{ \G^T: \hspace{-.7cm} &
 {\bullet^v} \ar@(u,ul)[] \ar@(d,dl)[] &{\, \bullet^{u } \,} \ar@/^/[l] \ar@/_/[l] &{\bullet^{u_{-1}}\,}\ar@/^/[l]^f \ar@/_/[l]_e  &{\bullet^{u_{-2}}\,}\ar@/^/[l] \ar@/_/[l] &{\dots}\ar@/^/[l] \ar@/_/[l]}
\end{align*}
}
\vspace{.2in}

\textup{In \cite[page 344]{DavR}, it is claimed that the projection P associated with the characteristic function of $\X \subseteq \X^T$ satisfies $P\O_{(\X^T, \sigma^T)}P= \O_{(\X, \sigma)}$ and so $\O_{(\X, \sigma)}$ is a corner of $\O_{(\X^T, \sigma^T)}$.
In our setting, this claim translates as follows: if $P=L_{u}+L_{v}$, then $P\O_{\G^{T}}P=\O_{\G}$. However this is not true. For instance, $P(L_f L_e^*)P=L_f L_e^* \notin \O_{\G}$.}
\end{example}

%%%%%%%%%%%%%%%%%%%%%%%%%%%%%%%%%%%%%
\subsubsection{Multivariable Dynamical Systems and Crossed Products by endomorphism}
%%%%%%%%%%%%%%%%%%%%%%%%%%%%%%%%%%%%%

We can describe the Cuntz-Pimsner algebra of an injective and non-degenerate multivariable system as a crossed product of a $\ca$-algebra $B$ by an endomorphism $\beta$. This idea appears first in \cite{DavR} for a different crossed product than the one presented here.
We start with the pertinent definitions.

%%%%%%%%%%%%%%%%%%%%%%%%%%%%%%%%%%%%%
\begin{definition} \label{covariant}
Let $B$ be a (not necessary unital) $\ca$-algebra and let $\beta$ be an injective endomorphism of $B$. A \textit{covariant} representation $(\pi, v)$ of the dynamical system $(B, \beta)$ consists of a non-degenerate $*$-representation $\pi$ of $B$ and an isometry $v$ satisfying
\begin{itemize}
\item[(i)] $\pi(\beta(b))=v\pi(b)v^*, \forall \, b \in B$, i.e. $v$ \textit{implements} $\beta$,
\item[(ii)] $v^*\pi(B)v\subseteq \pi(B)$, i.e. $v$ is \textit{normalizing} for $\pi(B)$,
\item[(iii)] $v^k(v^*)^k \pi(B) \subseteq \pi(B), \forall \, k \in \bbN.$
    \end{itemize}
    \end{definition}

The crossed product $B \times_{\beta} \bbN$ is the universal $\ca$-algebra associated with this concept of a covariant representation for $(B, \beta)$. Specifically, $B \times_{\beta} \bbN$ is generated by $B$ and $BV$, where $V$ is an isometry satisfying (i), (ii) and (iii) in Definition \ref{covariant} with $\pi=\id$.
Furthermore, for any covariant representation $(\pi, v)$ of $(B, \beta)$, there exists a $*$-homomorphism $\hat{\pi}\colon B \times_{\beta} \bbN \rightarrow \B(H)$ extending $\pi$ and satisfying $\hat{\pi}(bV) = \pi(b)v$, for all $b \in B$.

In the case where $B$ is unital, condition (iii) is redundant and this version of a crossed product by an endomorphism was introduced by Paschke~\cite{Pas}; in the generality presented here, it is new. It has the advantage that for any covariant representation of $(\pi, v)$ of $(B, \beta)$ admitting a gauge action, the fixed point algebra of $(\pi, v)$ equals $\pi(B)$ (see Remark \ref{fixed} below). This allows us to claim a gauge invariance uniqueness theorem for $B \times_{\beta} \bbN$: if $(\pi, v)$ is a faithful covariant representation of $(B, \beta)$ admitting a gauge action, then the $\ca$-algebra generated by $\pi(B)$ and $\pi(B)v$ is isomorphic to $B \times_{\beta} \bbN$.

%%%%%%%%%%%%%%%%%%%%%%%%%%%%%%%%%%%%%
\begin{remark}\label{fixed}
Observe that conditions (i), (ii) and (iii) in Definition~\ref{covariant} imply that $\hat{\pi}\left(B \times_{\beta} \bbN \right)$ is generated by polynomials of the form $\pi(b_0) + \sum_{k}\, \pi(b_k ) v^k +\sum_{l} \,(v^*)^l \pi( b_l)$,  \, $b_0 , b_k , b_l \in B$. Indeed, item (i) implies that $v\pi(b)=\pi(\beta(b))v$, thus $\pi(b)v^*=v^*\pi(\beta(b))$, for all $b\in B$ (since $B$ is self-adjoint). With this in hand and item (ii) of the definition we can see that the product $(v^*)^k\pi(b_k)\pi(b_l)v^l$ can be written as a trigonometric polynomial of the above form, for all $k,l\geq 0$. Item (iii) is used to show the same thing for the products $\pi(b_k)v^k(v^*)^l\pi(b_l)$, when $k\geq l$. Finally, for $k<l$ we use the fact that $v$ is an isometry and item (iii) so that
\begin{align*}
\pi(b_k)v^k(v^*)^l\pi(b_l) &= \pi(b_k)(v^*)^{l-k}v^l(v^*)^l\pi(b_l)\\
 &=\pi(b_k) (v^*)^{l-k} \pi(b') = (v^*)^{l-k}\pi(b''),
\end{align*}
which completes the argument.
\end{remark}

There is a related concept of a crossed product by an endomorphism which we now discuss. For a $\ca$-algebra $B$ and an injective endomorphism $\beta$, Stacey \cite{Stac} imposes on a covariant representation $(\pi, v)$ of $(B, \beta)$ only condition (i) from Definition~\ref{covariant}. He then defines the crossed product $B \rtimes_{\beta} \bbN$ to be the universal $\ca$-algebra associated with his concept of a covariant representation for $(B, \beta)$. Muhly and Solel have shown~\cite{MS2} that in the case where $B$ is unital, Stacey's crossed product is the Cuntz-Pimsner algebra of a certain correspondence. Using a gauge invariance uniqueness theorem one can prove that if the isometry V in $B \rtimes_{\beta} \bbN$ satisfies condition (ii) in Definition~\ref{covariant}, then $B \rtimes_{\beta} \bbN \simeq  B \times_{\beta} \bbN$.

%%%%%%%%%%%%%%%%%%%%%%%%%%%%%%%%%%%%%
\begin{theorem} \label{injectivemulti}
\cite[Theorem 4.6]{KakKats11} If $(A, \ga)$ is an injective multivariable system, then there exists a $\ca$-algebra $B$ and an injective endomorphism $\beta$ of $B$ so that $\O_{(A, \ga)}$ is isomorphic to the crossed product algebra $B\times_{\beta}\bbN$. Furthermore, if $A$ belongs to a class $\C$  which is invariant under direct limits and tensoring by $M_k(\bbC)$, $k\in \bbN$, then $B$ also belongs to $\C$.
\end{theorem}

Combining Theorem~\ref{injectivemulti} with Theorem~\ref{multitail} we obtain the following.

%%%%%%%%%%%%%%%%%%%%%%%%%%%%%%%%%%%%%
\begin{corollary} \label{main}
If $(A, \ga)$ is a multivariable system, then there exists a $\ca$-algebra $B$ and an injective endomorphism $\beta$ of $B$ so that $\O_{(A, \ga)}$ is isomorphic to a full corner of the crossed product algebra $B\times_{\beta}\bbN$. Furthermore, if $A$ belongs to a class $\C$  which is invariant under direct limits, quotients and tensoring by $M_k(\bbC)$, $k\in \bbN$, then $B$ also belongs to $\C$.
\end{corollary}

Finally, let us give a quick application of Theorem~\ref{injectivemulti}, that readily follows  from Paschke's result \cite{Pas} on the simplicity of $B \times_{\beta}\bbN$.

%%%%%%%%%%%%%%%%%%%%%%%%%%%%%%%%%%%%%
\begin{corollary}
Let $A$ be a UHF $\ca$-algebra and let $\ga=(\ga_1,\dots,\ga_n)$ be a multivariable system with $n\geq2$. If $\ga_i(1)=1$, for all $i=1,2,\dots,n$, then $\O_{(A,\ga)}$ is simple.
\end{corollary}

%%%%%%%%%%%%%%%%%%%%%%%%%%%%%%%%%%%%%
\subsection{Graph $\ca$-correspondences}
%%%%%%%%%%%%%%%%%%%%%%%%%%%%%%%%%%%%%

Graph $\ca$-correspondences can be used to summarize in a very delicate way properties and constructions of $\ca$-correspondences.

Let $\G=(\vrt,\edg,r,s)$ be a directed graph. Let $(A_p)_{p \in \vrt}$ be a family of $\ca$-algebras parameterized by the vertices of $\G$ and for each $e \in \edg$, we consider $X_e$ be a $A_{r(e)}$-$A_{s(e)}$-correspondences $X_e$.

Let $A_\G= c_0 (\, (A_p)_{p \in \vrt})$ denote the $c_0$-sum of the family $(A_p)_{p \in \vrt}$. Also, let $Y_0=c_{00}((X_e)_{e \in \edg})$ which is equipped with a $A_\G$-valued inner product
\[
\sca{ u, v  }(p)= \sum_{s(e)=p} \, \sca{u_e,v_e}_{A_p}, \quad p \in
\vrt,
\]
If $X_\G$ is the completion of $Y_0$ with respect to the inner product, then $X_\G$ is a $A_\G$-Hilbert module when equipped with the right action
\[
(ux)_e = u_ex_{s(e)}, \quad e \in \edg.
\]
It becomes a $\ca$-correspondence, when equipped with the $*$-homomorphism $\phi_{\G}: A_\G \rightarrow \L(X_\G)$, such that
\[
(\phi_{\G}(x)u)_e= \phi_e(x_{r(e)})(u_e), e\in \edg,
\]
for $x\in A_\G$ and $u \in X_\G$. The $\ca$-correspondence $X_\G$ over $A_\G$ is called \emph{the graph $\ca$-correspondence} with respect to the data
\[
\big\{\G, \{A_p\}_{p \in \vrt}, \{X_e\}_{e\in \edg} \big\}.
\]
In other words, a graph $\ca$-correspondence can be viewed as a graph $\G$ such that on every vertex sits a $\ca$-algebra, on each edge sits a $\ca$-correspondence and the actions and the inner product are defined via the information provided by the graph (they ``remember'' the form of the graph) as in the following image
\begin{align*}
\xymatrix{ & & & \ar@{-->}@/_/[dl]\\
& \bullet^{A_{r(e)}} \ar@{-->}@/_/[ul] \ar@{-->}@/_/[l] &
\bullet^{B_{s(e)}} \ar@/_/[l]_{X_e} & \ar@{-->}@/_/[l]}
\end{align*}

%%%%%%%%%%%%%%%%%%%%%%%%%%%%%%%%%%%%%
\begin{example}
Every $\ca$-correspondecne $X$ over $A$ can be visualized in the language of graph $\ca$-correspondences trivially as
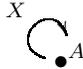
\begin{align*}
\xymatrix{ \bullet^A \ar@(l,u)[]^X }
\end{align*}
Moreover, every ${}_A X_B$ $\ca$-correpondence can be visualized as
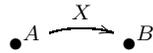
\begin{align*}
\xymatrix{ \bullet^A \ar@/^/[r]^X  & \bullet^{B}}
\end{align*}
\end{example}

The Fock space of a graph $\ca$-correspondence contains the \emph{path $\ca$-correspondences} $X_\mu$, for some path $\mu$. That is if $\mu=xe_n\dots e_1 y$ is a path of the graph $\G$ then $X_\mu:= X_{e_n} \otimes_{A_r(e_n)} \dots \otimes_{A_{r(e_2)}} X_{e_1}$ is a $A_x$-$A_y$-correspondence.

%%%%%%%%%%%%%%%%%%%%%%%%%%%%%%%%%%%%%
\begin{definition}
\textup{A graph $\ca$-correspondence is called \emph{commutative} if for any two paths $\mu= xe_n\dots e_1y$ and $\nu=xf_k\dots f_1y$, that have the same range and source, the corresponding \emph{path $\ca$-correspondences} $X_{e_n} \otimes \dots \otimes X_{e_1}$ and $X_{f_n}\otimes \dots X_{f_1}$ are unitarily equivalent.}
\end{definition}

For example if $E \approx R \otimes_B S$ for some $\ca$-correspondences ${}_A E_A, {}_A R_B$ and ${}_B S_A$, then one can form the commutative graph $\ca$-correspondence
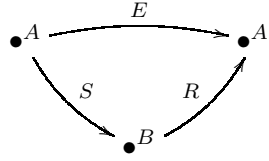
\begin{align*}
\xymatrix{ \bullet^A \ar@/^/[rr]^E \ar@/_/[dr]^S & & \bullet^A\\
& \bullet^B \ar@/_/[ur]^R &}
\end{align*}
Note that the graph $\ca$-correspondence that comes from this graph does not contain ${}_A E_A$ as a subcorrespondence, though. A graph $\ca$-correspondence that has additionally that property is the following
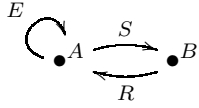
\begin{align*}
\xymatrix{\bullet^A \ar@(l,u)[]^E \ar@/^/[r]^S & \bullet^B
\ar@/^/[l]^R}
\end{align*}

%%%%%%%%%%%%%%%%%%%%%%%%%%%%%%%%%%%%%
\begin{example}
Let $\xtau$  be the correspondence that occurs from the adding-tail construction established in \cite{KakKats11}. It can be visualized in the language of graph $\ca$-correspondences as follows
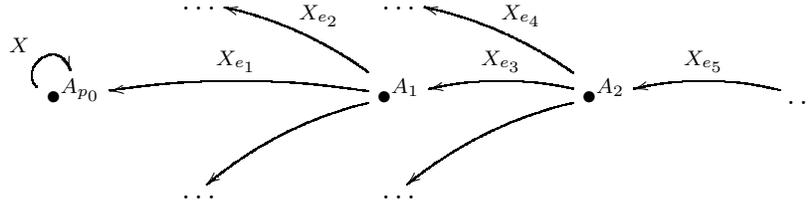
\begin{align*}
\xymatrix{& \dots & & \dots & & & & \\
\bullet^{A_{p_0}} \ar@(l,u)[]^{X} & & & \bullet^{A_1}
\ar@/_/[lll]_{X_{e_1}} \ar@/_/[ull]_{X_{e_2}} \ar@/_/[dll] & & \bullet^{A_2}
\ar@/_/[ll]_{X_{e_3}} \ar@/_/[ull]_{X_{e_4}} \ar@/_/[dll] & & \dots
\ar@/_/[ll]_{X_{e_5}}\\
& \dots & & \dots & & & &}
\end{align*}
where $A_{p_0} \equiv A$.
\end{example}

%%%%%%%%%%%%%%%%%%%%%%%%%%%%%%%%%%%%%
\begin{example}
\textup{(Muhly-Tomforde tail \cite{MuTom04}) When $X$ is visualized on a cycle graph, then Muhly-Tomforde tail produces the following graph $\ca$-correspondence}
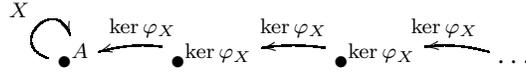
\begin{align*}
\xymatrix{ \bullet^A \ar@(l,u)[]^X & \bullet^{\ker\phi_X}
\ar@/_/[l]_{\ker\phi_X} & \bullet^{\ker\phi_X}
\ar@/_/[l]_{\ker\phi_X} & \cdots \ar@/_/[l]_{\ker\phi_X}}
\end{align*}
\end{example}

%%%%%%%%%%%%%%%%%%%%%%%%%%%%%%%%%%%%%
\begin{example}
If $(X_{\ga}, A)$ is the $\ca$-correspondence canonically associated with a dynamical system $(A, \ga)$, then the tail produced by Muhly and Tomforde for $(X_{\ga}, A)$ comes from the graph
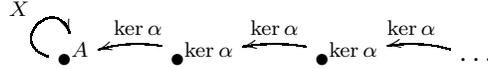
\begin{align*}
\xymatrix{ \bullet^A \ar@(l,u)[]^X & \bullet^{\ker\ga}
\ar@/_/[l]_{\ker\ga} & \bullet^{\ker\ga}
\ar@/_/[l]_{\ker\ga} & \cdots \ar@/_/[l]_{\ker\ga}}
\end{align*}
\end{example}

%%%%%%%%%%%%%%%%%%%%%%%%%%%%%%%%%%%%%
\begin{example}
If $(X_{\ga}, A)$ is the $\ca$-correspondence canonically associated with a dynamical system $(A, \ga)$, then the appropriate tail for $(X_{\ga}, A)$ comes from the graph
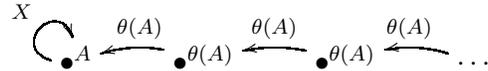
\begin{align*}
\xymatrix{ \bullet^A \ar@(l,u)[]^X & \bullet^{\theta(A)}
\ar@/_/[l]_{\theta(A)} & \bullet^{\theta(A)}
\ar@/_/[l]_{\theta(A)} & \cdots \ar@/_/[l]_{\theta(A)}}
\end{align*}
where $\theta \colon A \longrightarrow M(\ker \ga)$ is the map that extends the natural inclusion $\ker \ga \subseteq M(\ker \ga))$ in the multiplier algebra. Therefore, the $\ca$- correspondence $(Y,\phi_Y, B)$ is the $\ca$-correspondence canonically associated with a dynamical system $(B, \beta)$, where $B=A\oplus c_0(\theta(A))$ and $\beta(a,(x_n))=(\alpha(a),\theta(a),(x_n))$.
\end{example}

%%%%%%%%%%%%%%%%%%%%%%%%%%%%%%%%%%%%%
\begin{example}
If $(X_{\ga}, A)$ is the $\ca$-correspondence canonically associated with a $\ca$-dynamical system $(A, \ga_1, \ga_2)$, then the appropriate tail for $(X_{\ga}, A)$ comes from the graph
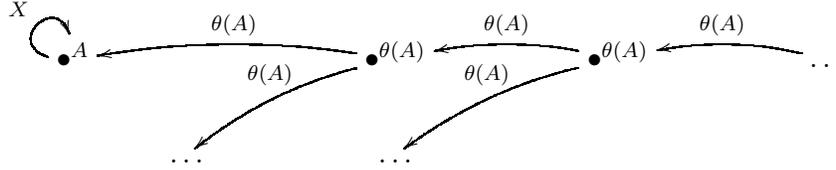
\begin{align*}
\xymatrix{\bullet^{A} \ar@(l,u)[]^{X} & & & \bullet^{\theta(A)}
\ar@/_/[lll]_{\theta(A)} \ar@/_/[dll]_{\theta(A)} & & \bullet^{\theta(A)}
\ar@/_/[ll]_{\theta(A)}  \ar@/_/[dll]_{\theta(A)} & & \dots
\ar@/_/[ll]_{\theta(A)}\\
& \dots & & \dots & & & &}
\end{align*}
where $\J= \ker \ga_1 \cap \ker\ga_2$ and $\theta \colon A \longrightarrow M(\J)$ is the map that extends the natural inclusion $\J \subseteq M(\J)$ in the multiplier algebra.
\end{example}

%%%%%%%%%%%%%%%%%%%%%%%%%%%%%%%%%%%%%
\section{Shift Equivalence Problem}
%%%%%%%%%%%%%%%%%%%%%%%%%%%%%%%%%%%%%

In his pioneering paper \cite{Wil74} Williams studied certain notions of relations for the class of matrices with non-negative integer entries. We say that two such matrices $E$ and $F$ are \emph{elementary strong shift equivalent}, and write $E\sse F$, if there are two matrices $R$ and $S$ such that $E = RS$ and $F=SR$. Whereas this relation is symmetric, it may not be transitive. In \cite[Example 2]{Wil74} Williams gives the following counterexample; let
\begin{align*}
\left[\begin{array}{cc} 10 & 2\\ 2 & 1
\end{array}\right] \sse \left[\begin{array}{cc} 5 & 6\\ 4 & 6
\end{array}\right] \sse \left[\begin{array}{cc} 9 & 4\\ 3 & 2
\end{array}\right],
\end{align*}
but $\left[\begin{array}{cc} 10 & 2\\ 2 & 1 \end{array}\right]$ is not elementary strong shift equivalent to $\left[\begin{array}{cc} 9 & 4\\ 3 & 2 \end{array}\right]$. The transitive closure $\sse$, denoted by $\tsse$, implies that $E \tsse F$ if there is a sequence of matrices $T_i$, $i=0,\dots, n$, such that $E=T_0$, $F=T_n$ and $T_i \tsse T_{i+1}$. Williams also defined a third relation, which is proved to be transitive. We say that $E$ is \emph{shift equivalent} to $F$, and write $E \se F$, if there are $R,S$ such that $E^n = RS$, $F^n=SR$ and $ER=SF$, $FR=SE$ for some $n\in \bbN$.

A purpose of \cite{Wil74} was to prove that the relations $\tsse$ and $\se$ are equivalent. Unfortunately, an error in \cite{Wil74} made invalid the proof of a key lemma, and this task remained unsolved for over than 20 years, known as \emph{Williams' Conjecture}. The research interest in this area contributed to the growth of symbolic dynamics and to the search of (complete) invariants for both the equivalence relations. A major change was made by Kim and Roush in \cite{KimRoush99}, where they proved that Williams' Conjecture was false for the class of non-negative integral matrices. Their work suggests that Williams' Conjecture can be renamed as \emph{Shift Equivalence Problem}; i.e. is $\tsse$ equivalent to $\se$ for a class $\fS$. This formulation is a little vague as one has to extend the definition of the relations described above to \emph{a} class $\fS$.\\

The notion of elementary and strong shift equivalence for $\ca$- correspondences was studied by Muhly, Tomforde and Pask \cite{MuPasTom08}. In addition they prove that strong shift equivalence of $\ca$-correspondences implies the Morita equivalence of the associated Cuntz-Pimsner algebras, thus extending classical results of Cuntz and Krieger \cite{CunKri80}, Bates \cite{Bat02} and Drinen and Sieben \cite{DriSie01} for graph $\ca$-algebras.
The concept of shift equivalence has been studied extensively from both the dynamical and the ring theoretic viewpoint. (See \cite{Wagoner} for  a comprehensive exposition.) In general, shift equivalence  has been  recognized to be a more manageable  invariant than strong shift equivalence, as it is decidable over certain rings \cite{KimRoushdecide}. Unlike  strong shift equivalence, the study of shift equivalence, from the viewpoint of $\ca$-correspondences, has  been met with limited success \cite{Matsumoto}. (Other operator theoretic viewpoints however have been quite successful \cite{Krieger}.)
The concept of strong Morita equivalence for $\ca$-correspondences was first developed and studied by Abadie, Eilers and Exel \cite{AbaEilEx98} and Muhly and Solel \cite{MuhSol00}, and plays the role of a generalized Conjugacy (see Example \ref{conjug} below). Among others these authors show that if two $\ca$-correspondences are strong Morita equivalent then the associated Cuntz-Pimsner algebras $\O_E$ and $\O_F$ are (strong) Morita equivalent as well. Let us give the definitions.

%%%%%%%%%%%%%%%%%%%%%%%%%%%%%%%%%%%%%
\begin{definition}
\textup{Let the $\ca$-correspondences ${}_A E_A$ and ${}_B F_B$. Then we say that}
\begin{enumerate}
\item \textup{\emph{$E$ is Morita equivalent to $F$}, and we write $E\sme F$, if there is an imprimitivity bimodule ${}_A M_B$ such that $E \otimes_A M = M\otimes_B F$.}
\item \textup{\emph{$E$ is elementary strong shift equivalent to
$F$}, and we write $E\sse F$, if there are ${}_A R_B$ and ${}_B S_A$ such that $E= R\otimes_A
S$ and $F=S \otimes_B R$,}
\item \textup{\emph{$E$ is strong shift equivalent to $F$}, and we write $E \tsse F$, if there are $T_i$, $i=0,\dots n$, such that $T_0=E, T_n=F$ and $T_i \sse T_{i+1}$,}
\item \textup{\emph{$E$ is shift equivalent to $F$ with lag $m$}, and we write $E \se F$, if there are ${}_A R_B$ and ${}_B S_B$ such that $E^{\otimes m}= R\otimes_B S$, $F^{\otimes m}= S \otimes_A R$ and $E\otimes_A R=R \otimes_B F$, $S\otimes_A E= F\otimes_B S$.}
\end{enumerate}
\end{definition}

%%%%%%%%%%%%%%%%%%%%%%%%%%%%%%%%
\begin{example}\label{conjug}
\textup{Let ${}_A E_A$ and ${}_B F_B$ be $\ca$-correspondences arising from two dynamical systems $(A,\ga)$ and $(B,\beta)$ and assume that $E \sme F$ via an imprimitivity bimodule $M$. If we wish $M$ to arise in a similar way then it should be the $\ca$-correspondence associated to a $*$-isomorphism $\gamma\colon A \rightarrow B$. As showed in Example \ref{Ex:ten dyn}, $E\otimes_A M$ is the $\ca$-correspondence associated to $\gamma\circ \ga\colon A \rightarrow B$ and $M\otimes_B B$ is the $\ca$-correspondence associated to $\beta\circ \gamma\colon A \rightarrow B$. Therefore the unitary equivalence $E \otimes_A M \approx M\otimes_B B$ induces a unitary $u\in \L(B)=B$, such that $a\cdot u(1_B)= u(a\cdot 1_B)$, for all $a\in A$; equivalently that $\beta\circ \gamma(a)u=u \gamma\circ \ga(a)$, for all $a\in A$, hence that the systems $(A,\ga)$ and $(B,\beta)$ are (outer) conjugate.}
\end{example}

In \cite{KakKats12} we studied these relations and the interaction between them. Note that $\sme, \tsse$ and $\se$ are equivalence relations for non-degenerate $\ca$- correspondences. One of our main result is the following.

%%%%%%%%%%%%%%%%%%%%%%%%%%%%%%%%%%%%%
\begin{theorem}
\cite{KakKats12} The Shift Equivalence Problem Problem is true for the class of imprimitivity bimodules. In particular, the relations $\sme, \sse, \tsse, \se$ coincide.
\end{theorem}

A weaker version of the previous Theorem holds in general.

\begin{theorem}
\cite{KakKats12} Let the $\ca$-correspondences ${}_A E_A$ and ${}_B F_B$; then
\begin{align*}
\xymatrix{ E \sme F \ar@{=>}[r] & E \sse F  \ar@{=>}[r] & E\tsse F \ar@{=>}[r] &  E \se F}.
\end{align*}
\end{theorem}

In \cite{MuPasTom08} Muhly, Pask and Tomforde have provided a number of counterexamples showing that Morita equivalence of $\ca$-correspondences differs from the elementary strong shift equivalence. The previous result shows that it is in fact stronger.\\

One of the basic tools we use in \cite{KakKats12} is the \emph{Pimsner dilation} of an injective $\ca$-correspondence $X$ to a Hilbert bimodule $X_\infty$. This construction was first introduced by Pimsner in \cite{Pim}. In \cite[Appendix A]{KakKats11} we revisited this construction by using direct limits. Consider the directed system
\[
A \stackrel{\rho_0}{\longrightarrow} \L(X)
\stackrel{\rho_1}{\longrightarrow}  \L(X^{\otimes 2})
\stackrel{\rho_2} \longrightarrow \cdots,
\]
where
\begin{align*}
&\rho_0=\phi_X \colon A=\L(A)\longrightarrow \L(X),\\
&\rho_n \colon \L(X^{\otimes n}) \longrightarrow \L(X^{\otimes
n+1})\colon r \longmapsto r\otimes \id_X , \, n\geq 1,
\end{align*}
and let $ A_\infty$ be the $\ca$-subalgebra of $B=\varinjlim (\L(X^{\otimes n}), \rho_n)$ that is generated by the copies of $\K(X^{\otimes n})$, for $n\in \bbZ_+$. Consider also the directed system of Banach spaces
\[
\L(A,X) \stackrel{\sigma_0}{\longrightarrow} \L(X,X^{\otimes 2})
\stackrel{\sigma_1}{\longrightarrow} \cdots ,
\]
where
\begin{align*}
\sigma_n \colon \L(X^{\otimes n}, X^{\otimes n+1}) \rightarrow
\L(X^{\otimes n+1}, X^{\otimes n+2}) \colon  s \mapsto s\otimes
\id_X , \, n\geq 1,
\end{align*}
and let $X_\infty$ be the Banach subalgebra of $Y= \varinjlim (\L(X^{\otimes n}, X^{\otimes n+1}), \gs_n)$ generated by the copies of $\K(X^{\otimes n}, X^{\otimes n+1})$, for $n\in \bbZ_+$. Note that the map
\[
\partial \colon X \longrightarrow
\L(A,X) \colon \xi \longmapsto \partial_{\xi},
\]
where $\partial_{\xi}(a)=\xi a$, $\xi \in X$, maps a copy of $X$ isometrically into $\K(A,X) \subseteq X_\infty$. In particular, if $\phi_X(A) \subseteq \K(X)$, then one can verify that $X_\infty = \varinjlim (\K(X^{\otimes n}, X^{\otimes n+1}), \gs_n)$ and $A_\infty = \varinjlim (\K(X^{\otimes n}), \rho_n)$. Thus, in this case, $X_\infty$ is a full left Hilbert bimodule.

If $r\in \L(X^{\otimes n})$, $s\in \L(X^{\otimes n},X^{\otimes n+1})$ and $[r], [s]$ are their equivalence classes in $B$ and $Y$ respectively, then we define $[s]\cdot [r]:= [sr]$. From this, it is easy to define a right $B$-action on $Y$. Similarly, we may define a $B$-valued right inner product on $Y$ by setting
\begin{align*}
\sca{[s'],[s]}_Y \equiv [(s')^*s] \in B.
\end{align*}
for $s, s' \in \L(X^{\otimes n},X^{\otimes n+1})$, $n \in \bbN$, and then extending to $Y \times Y$. Finally we define a $*$-homomorphism $\phi_Y \colon B \rightarrow \L(Y)$ by setting
\[
\phantom{XXX} \phi_Y([r])([s]) \equiv [rs], \quad r\in \L(X^{\otimes
n}), s\in \L(X^{\otimes n-1},X^{\otimes n}), n\geq0
\]
and extending to all of $B$ by continuity. We therefore have a left $B$-action on $Y$ and thus $Y$ becomes a $\ca$-correspondence over $B$.

The following diagrams depict the above construction in a heuristic way: the right action is ``defined'' through the diagram
\begin{align*}
\xymatrix{ A \ar[r]^{\phi_X\equiv \rho_0} \ar@{-->}[d]_{\cdot} &
\L(X) \ar[r]^{\rho_2} \ar@{-->}[d]_{\cdot} & \L(X^{\otimes 2})
\ar[r]^{\rho_3} \ar@{-->}[d]_{\cdot} &
\dots \ar[r] \ar@{-->}[d]_{\cdot} & B \ar@{-->}[d]_{\cdot} \\
\L(A,X) \ar[r]^{\sigma_0}  & \L(X,X^{\otimes 2}) \ar[r]^{\sigma_1} &
\L(X^{\otimes 2},X^{\otimes 3}) \ar[r]^{\phantom{XXX} \sigma_3} &
\dots \ar[r] & Y}
\end{align*}
while the left action is ``defined'' through the diagram
\begin{align*}
 \xymatrix{A \ar[r]^{\phi_X\equiv \rho_0}
 & \L(X) \ar[r]^{\rho_2} \ar@{-->}[dl]_{\id \cdot} &
 \L(X^{\otimes 2}) \ar[r]^{\rho_3} \ar@{-->}[dl]_{\id \cdot} &
 \dots \ar@{-->}[dl]_{\id \cdot} \ar[r] & B \ar@{-->}[d]_{\phi_{Y}}  \\
 \L(A,X) \ar[r]^{\sigma_0}  & \L(X,X^{\otimes 2})
 \ar[r]^{\sigma_1} & \L(X^{\otimes 2},X^{\otimes 3}) \ar[r]^{} &
 \dots \ar[r] & Y}
\end{align*}

For a proof of the following Theorem see \cite[Theorem 2.5]{Pim} or \cite[Theorem 6.6]{KakKats11}. The main difference between the two approaches is that in \cite{KakKats11} we have represented $X \otimes_A A_\infty$ as the direct limit $X_\infty$, thus having no concern in checking the form of the tensor product. Moreover, it appears that ``tensoring'' $X$ with $A_\infty$ is equivalent to multiplying $X$ with $A_\infty$ in $\varinjlim (\L(X^{\otimes n}, X^{\otimes n+1}), \gs_n)$.

%%%%%%%%%%%%%%%%%%%%%%%%%%%%%%%%%%%%%
\begin{theorem}\label{pim}
\cite[Theorem 2.5]{Pim} \cite[Theorem 6.6]{KakKats11} Let ${}_A X_A$ be an injective $\ca$-correspondence and let $X_{\infty}$ be the $A_{\infty}$-correspondence defined above. Then $X_{\infty}$ is an essential Hilbert bimodule and its Cuntz-Pimsner algebra $\O_{X_\infty}$ is $*$-isomorphic to $\O_X$.
\end{theorem}

The idea of using direct limits for Pimsner dilation is in complete analogy to the direct limit process for dynamical systems (see \cite{kakarkats10}).

\begin{example}\label{Ex:dil dyn}
\textup{Let $(A,\ga)$ denote a dynamical system where $\ga$ is a $*$- injective endomorphism of $A$. We can define the direct limit dynamical system $(A_\infty,\ga_\infty)$ by}
\begin{align*}
 \xymatrix{
  A \ar[r]^{\ga} \ar[d]^\ga &
  A \ar[r]^{\ga} \ar[d]^\ga &
  A \ar[r]^{\ga} \ar[d]^\ga &
  \cdots \ar[r] &
  A_\infty \ar[d]^{\ga_\infty} \\
  A \ar[r]^\ga &
  A \ar[r]^\ga &
  A \ar[r]^\ga &
  \cdots \ar[r] &
  A_\infty
 }
\end{align*}
\textup{The limit map $\ga_\infty$ is an automorphism of $A_\infty$ and extends $\ga$ (note that $A$ imbeds in $A_\infty$ since $\ga$ is injective). Then the $A_\infty$-$A_\infty$-correspondence $X_{\ga_\infty}$, is the Pimsner dilation of $X_\ga$.}
\end{example}

The main question imposed in \cite{KakKats12} was the following: is it true that $E_\infty \sim F_\infty$ (where $\sim$ may be $\sme, \sse, \tsse, \se$) if $E\sim F$ (in the same way)? When $E$ and $F$ are non-degenerate and regular $\ca$-correspondences we get the following result, by making use of the generalized notion of dilations.

%%%%%%%%%%%%%%%%%%%%%%%%%%%%%%%%%%%%%
\begin{theorem}\cite{KakKats12}
Let ${}_A E_A$ and ${}_B F_B$ be non-degenerate and regular $\ca$-correspondences. If $E\sim F$ $($where $\sim$ may be $\sme, \sse, \tsse, \se$), then $E_\infty \sim F_\infty$ (in the same way$)$.
\end{theorem}

For full right, non-degenerate, regular $\ca$-correspondences we have the following.

%%%%%%%%%%%%%%%%%%%%%%%%%%%%%%%%
\begin{theorem}\label{full}
\cite{KakKats12} Let ${}_A E_A$ and ${}_B F_B$ be full right, non-degenerate, regular $\ca$-correspondences. Then the above scheme holds
\begin{align*}
\xymatrix{
 E\sme F     \ar@{=>}[d] \ar@{=>}[r]
 & E \sse F  \ar@{=>}[d] \ar@{=>}[r]
 & E\tsse F  \ar@{=>}[d] \ar@{=>}[r]
 & E\se F    \ar@{=>}[d] \\
 E_\infty \sme F_\infty     \ar@{<=>}[r]
 & E_\infty \sse F_\infty  \ar@{<=>}[r]
 & E_\infty \tsse F_\infty  \ar@{<=>}[r]
 & E_\infty \se F_\infty }
\end{align*}
\end{theorem}

The vertical arrows in Theorem \ref{full} are not equivalences in general. Indeed, if $E_\infty \sim F_\infty$ implied $E\sim F$ then, in particular $E \se E_\infty$. The following Theorem shows that this happens only trivially, i.e. when $E=E_\infty$.

%%%%%%%%%%%%%%%%%%%%%%%%%%%%%%%%
\begin{theorem}\label{T:se not inv}
\cite{KakKats12} Let ${}_A E_A$ be a full, non-degenerate and regular $\ca$-correspondence. If  $E \se E_\infty$ then $E$ is an imprimitivity bimodule.
\end{theorem}

It would be interesting if we could prove the validity of the Shift Equivalence Problem for this class of $\ca$-correspondences (as, after imprimitivity bimodules, it is the next best thing). An obstacle that prevents the construction of counterexamples that would give a negative answer to the Shift Equivalence Problem for non-degenerate and regular correspondences is that the theory of invariants of correspondences is rather poor. The best results (to our opinion) obtained so far are those appearing in \cite{MuhSol00, MuPasTom08} which we state for sake of completeness. Note that the additional item (3) below is an immediate consequence of item (2) and item (4) below is an immediate consequence of Theorem \ref{pim} and Theorem \ref{full}.

%%%%%%%%%%%%%%%%%%%%%%%%%%%%%%%%%%%%%
\begin{theorem}\label{mor equiv}
\textup{(1)} \textup{\cite[Theorem 3.2]{MuhSol00}} Let ${}_A E_A$ and ${}_B F_B$ be non- degenerate, injective $\ca$-correspondences. If $E \sme F$ then the corresponding Toeplitz-Cuntz-Pimsner algebras and Cuntz-Pimsner algebras are strong Morita equivalent as $\ca$-algebras, and the corresponding tensor algebras are strong Morita equivalent in the sense of \cite{BleMuPau00}.

\textup{(2)} \textup{\cite[Theorem 3.14]{MuPasTom08}} Let ${}_A E_A$ and ${}_B F_B$ be non-degenerate, regular $\ca$-correspondences. If $E\sse F$, then $\O_E \sme \O_F$.

\textup{(3)} \cite{KakKats12} Let ${}_A E_A$ and ${}_B F_B$ be non-degenerate, regular $\ca$- correspondences. If $E\tsse F$, then $\O_E \sme \O_F$.

\textup{(4)} \cite{KakKats12} Let ${}_A E_A$ and ${}_B F_B$ be full right, non-degenerate, regular $\ca$-correspondences. If $E\se F$, then $\O_E \sme \O_F$.
\end{theorem}

In particular we obtain the following result for Cuntz-Krieger $\ca$-algebras.

%%%%%%%%%%%%%%%%%%%%%%%%%%%%%%%%%%%%%
\begin{corollary} \label{graph}
Let  $\G$ and $\G'$ be finite graphs with no sinks or sources and let $A_{\G}$ and $A_{\G'}$ be their adjacent matrices. If $A_{\G} \se A_{\G'}$, in the sense of Williams, then the Cuntz-Krieger $\ca$-algebras $\O_{\G}$ and $\O_{\G'}$ are strong Morita equivalent.
\end{corollary}

There is also a direct application to unital injective dynamical systems.

%%%%%%%%%%%%%%%%%%%%%%%%%%%%%%%%%%%%%
\begin{corollary} \label{dyn sys}
Let $(A,\ga)$ and $(B,\beta)$ be unital injective dynamical systems. If $X_\ga \se X_\beta$, then $X_{\ga_\infty} \se Y_{\beta_\infty}$ and the crossed products $A_\infty \rtimes_{\ga_\infty} \bbZ$ and $B\rtimes_{\beta_\infty} \bbZ$ are strong Morita equivalent.
\end{corollary}

Theorem \ref{mor equiv} shows that Cuntz-Pimsner algebras is a rather coarse invariant. After all, Cuntz-Pimsner algebras are not a complete invariant for the restriction of the relations $\sme$, $\sse$, $\tsse$ or $\se$ to subclasses of $\ca$-correspondences. For example, let $(A,\ga)$ be the dynamical system constructed by Hoare and Parry in \cite{HoaPar66}. Then $\ga$ is a $*$-isomorphism and $\ga$ is not conjugate to its inverse. If $E$ is the $\ca$-correspondence of $(A,\ga)$ and $F$ is the $\ca$-correspondence of $(A,\ga^{-1})$, then there is not a $\ca$-correspondence $M$ of a dynamical system $(B,\beta)$ such that $E\otimes M \approx M \otimes F$, because then the two dynamical systems would be conjugate. But $\O_A= A\rtimes_\ga \bbZ$ is always $*$-isomorphic to $\O_B=A \rtimes_{\ga^{-1}} \bbZ$.

On the other hand the tensor algebras of $\ca$-correspondences may be more eligible. They were used in \cite{MuPasTom08} to show that $\sse$ does not imply $\sme$ and they provide a complete invariant for the conjugacy problem for dynamical systems in various cases, as shown by Davidson and Katsoulis \cite{DK, DKsimple}, and recently by Davidson and Kakariadis \cite{DavKak12}. Moreover, for aperiodic $\ca$-correspondences, Muhly and Solel \cite{MuhSol00} prove that $\sme$ is equivalent to strong Morita equivalence of the tensor algebras in the sense of \cite{BleMuPau00}. Along with the further investigation of tensor algebras, it is natural to suggest to work on the development of other invariants such as periodicity, existence of cycles, saturated and/or hereditary submodules etc.

%%%%%%%%%%%%%%%%%%%%%%%%%%%%%%%%%%%%%
\section{Hilbert Bimodules}\label{hilbert bimodules}
%%%%%%%%%%%%%%%%%%%%%%%%%%%%%%%%%%%%%

Recall that if $X$ is a Hilbert bimodule then
\begin{align*}
 I_X=\overline{\Span}\{\lsca{\xi,\eta}: \xi,\eta \in X\}.
\end{align*}
If a Hilbert $A$-bimodule $X$ is considered as a $\ca$-correspondence over $A$, then $J_X=I_X$. Hence, for $\xi,\eta\in X$, the element $\lsca{\xi,\eta}\in A$ is identified with the unique element $a\in J_X$ such that $\phi_X(a)=\Theta^X_{\xi,\eta}$. The converse is also true.

The following result is well-known.

%%%%%%%%%%%%%%%%%%%%%%%%%%%%%%%%%%%%%
\begin{proposition}\label{equivalence via compacts}
Let $X$ be a $\ca$-correspondence over $A$. Then the following are equivalent
\begin{enumerate}
\item $X$ is a bimodule,
\item $K(X) \subseteq \phi_X(A)$ and
$\phi_X^{-1}(\K(X))=\ker\phi_X \oplus J_X$,
\item the restriction of $\phi_X$ to $J_X$ is a
$*$-isomorphism onto $\K(X)$.
\end{enumerate}
In particular, if $\phi_X$ is injective then $X$ is a Hilbert bimodule if and only if $\K(X) \subseteq \phi_X(A)$.
\end{proposition}

It turns out that the property of a $\ca$-correspondence being a Hilbert bimodule has an important non-selfadjoint operator algebra manifestation. We remind that an operator algebra $\fA$ is called \emph{Dirichlet} if $\fA+\fA^*$ is dense (via a completely isometric homomorphism) in $\cenv(\fA)$.

%%%%%%%%%%%%%%%%%%%%%%%%%%%%%%%%%%%%%
\begin{theorem}
\cite{KakKats12-2} Let $X$ be a $\ca$-correspondence over $A$. Then the following are equivalent:
\begin{enumerate}
\item $X$ is a Hilbert bimodule,
\item $\psi_t(\K(X)) \subseteq \pi(A)$, for any injective covariant representation
$(\pi,t)$ that admits a gauge action,
\item the tensor algebra $\T_X^+$ has the Dirichlet property.
\end{enumerate}
\end{theorem}

The above result allows us to correct a misconception regarding semicrossed products and Dirichlet algebras.

%%%%%%%%%%%%%%%%%%%%%%%%%%%%%%%%%%%%%
\begin{corollary} Let $(A, \ga)$ be a dynamical system. Then the semicrossed product $A\times_{\ga} \bbZ_+$ has the Dirichlet property if and only if $\ga$ is surjective and $\ker\ga$ is orthocomplemented in $A$.
\end{corollary}

In particular, when $\ga$ is injective we deduce that the semi-crossed product has the Dirichlet property if and only if $\ga$ is onto (thus a $*$-isomorphism). Thus \cite[Proposition 3]{Dun08} is false. Nevertheless the main results of \cite{Dun08} are correct since they do not require that Proposition. It remains of interest though to determine whether a semi-crossed product has the unique extension property, since this would allow us to extend Duncan's results to non-commutative dynamics.

Recall that an operator algebra $\fA$ is said to have the \emph{the unique extension property} if the restriction of every faithful representation of $\cenv(\fA)$ to $\fA$ is maximal, i.e. it has no non-trivial dilations.

%%%%%%%%%%%%%%%%%%%%%%%%%%%%%%%%%%%%%
\begin{theorem}
\cite{KakKats12-2} Let $(A,\ga)$ be a unital injective dynamical system of a $\ca$-algebra. Then the semicrossed product $A\times_{\ga} \bbZ_+$ has the unique extension property.
\end{theorem}

The unique extension property of an operator algebra (or in general an operator space) implies the existence of the Choquet boundary in the sense of Arveson \cite{Ar08}, i.e. the existence of sufficiently many irreducible representations such that their restriction is maximal. Indeed, let $PS(\cenv(\fA))$ be the set of the pure states of $\cenv(\fA)$, and let $\Pi=\oplus_{\tau \in PS(\cenv(\fA))} \pi_\tau$ be the free atomic representation of $\cenv(\fA)$. Then $\Pi$ is faithful on $\cenv(\fA)$, hence by the unique extension property its restriction to $\fA$ is maximal. Moreover, every $\pi_\tau$ is maximal as a direct summand of a maximal representation. Hence,
\begin{align*}
\nor{[x_{ij}]}=\nor{[\Pi(x_{ij})]}=\sup\{\nor{[\pi_\tau(x_{ij})]}: \tau \in PS(\cenv(\fA))\},
\end{align*}
for all $[x_{ij}] \in M_\nu(\fA)$ and $\nu\in \bbN$.

The existence of the Choquet boundary for separable operator systems (or operator spaces) was proved by Arveson in \cite{Ar08} and it is still an open problem for the non-separable cases. Recently it was proved by Kleski in \cite{Kle11} that the supremum above can be replaced by a maximum, at least for the separable cases, where Arveson's Theorem applies. Note that the semicrossed products can give examples of non-separable operator algebras that have a Choquet boundary.

Finally, we have a result that relates our adding of a tail process to the concept of a Hilbert bimodule.

%%%%%%%%%%%%%%%%%%%%%%%%%%%%%%%%%%%%%
\begin{theorem}\label{dirichlet bimodules}
\cite{KakKats12-2} Let $X$ be a non-injective $\ca$-correspondence. Then the graph $\ca$-correspondence $X_\tau$, as defined in section \ref{adding tails}, is an (essential) Hilbert bimodule if and only if $X$ is a Hilbert bimodule and $|s^{-1}(p)|=|r^{-1}(p)|=1$ for every $p \neq p_0$.
\end{theorem}

%%%%%%%%%%%%%%%%%%%%%%%%%%%%%%%%%%%%%

\end{document}